\documentclass[11pt]{article} 
\usepackage{amssymb,amsmath,mathtools}
\usepackage{color}
\usepackage[OT2, T1]{fontenc} 
\usepackage[utf8]{inputenc} 

\usepackage[noadjust]{cite} 

\usepackage[colorlinks=true, urlcolor=green!45!black, citecolor=red, linkcolor=blue,
pdfpagelabels, bookmarksnumbered, bookmarksopen]{hyperref} 

\usepackage[hyperpageref]{backref} 

\newcommand{\nequiv}{\not \equiv}

\usepackage[italian, english]{babel} 

\allowdisplaybreaks 

\def\claim#1.{\noindent {\bf #1.}}

\def\flushright#1{{\unskip\nobreak\hfil\penalty50\hskip2em\hbox{}\nobreak\hfil%
#1\parfillskip=0pt\finalhyphendemerits=0\par}}

\def\bull{\vrule height 1.8ex width 1.0ex depth .1ex }
\def\QED{\ifmmode\eqno\hbox{$\bull$}\else\flushright{\hbox{$\bull$}}\fi}

\newcommand{\parag}[1]{\left\{ \begin{aligned} #1 \end{aligned}\right.} 

\usepackage{scalerel}
\usepackage{tikz}
\usetikzlibrary{svg.path}
\definecolor{orcidlogocol}{HTML}{A6CE39}
\tikzset{
 orcidlogo/.pic={
 \fill[orcidlogocol] svg{M256,128c0,70.7-57.3,128-128,128C57.3,256,0,198.7,0,128C0,57.3,57.3,0,128,0C198.7,0,256,57.3,256,128z};
 \fill[white] svg{M86.3,186.2H70.9V79.1h15.4v48.4V186.2z}
 svg{M108.9,79.1h41.6c39.6,0,57,28.3,57,53.6c0,27.5-21.5,53.6-56.8,53.6h-41.8V79.1z M124.3,172.4h24.5c34.9,0,42.9-26.5,42.9-39.7c0-21.5-13.7-39.7-43.7-39.7h-23.7V172.4z}
 svg{M88.7,56.8c0,5.5-4.5,10.1-10.1,10.1c-5.6,0-10.1-4.6-10.1-10.1c0-5.6,4.5-10.1,10.1-10.1C84.2,46.7,88.7,51.3,88.7,56.8z};
 }
}
\newcommand\orcidicon[1]{\href{https://orcid.org/#1}{\mbox{\scalerel*{
\begin{tikzpicture}[yscale=-1,transform shape]
\pic{orcidlogo};
\end{tikzpicture}
}{|}}}}

\newtheorem{Theorem}{Theorem}[section]
\newtheorem{Proposition}[Theorem]{Proposition}
\newtheorem{Corollary}[Theorem]{Corollary}
\newtheorem{Lemma}[Theorem]{Lemma}
\newtheorem{Remark}[Theorem]{Remark}

\newcommand{\R}{\mathbb{R}}
\newcommand{\N}{\mathbb{N}}
\newcommand{\Z}{\mathbb{Z}}

\newcommand{\mc}[1]{\mathcal{#1}}
\newcommand{\wto}{\rightharpoonup}

\def\half{{1\over 2}}

\def\abs#1{|#1|}
\def\pabs#1{\left|{#1}\right|}
\def\norm#1{\|#1\|}

\def\epsilon{\varepsilon}
\def\eps{\varepsilon}

\newcommand{\G}{{\mathbb G}}
\newcommand{\sgn}{{\rm sgn}}
\newcommand{\dista}{{\rm dist}}

\def\supp{\mathop{\rm supp}}

\def\RRint{\iint_{\R^N\times\R^N}}
\def\meas{\mathop{\rm meas}\nolimits} 
\def\calD{{\cal D}}

\newcommand{\dive}{{\rm div}}

\begin{document}



\title
{Infinitely many free or prescribed mass solutions for fractional Hartree equations and Pohozaev identities}

\author{
		Silvia Cingolani 
		\orcidicon{0000-0002-3680-9106}
		\\ \normalsize{Dipartimento di Matematica}
		\\ \normalsize{Universit\`{a} degli Studi di Bari Aldo Moro}
		\\ \normalsize{Via E. Orabona 4, 70125 Bari, Italy}
		\\ \normalsize{\href{mailto:silvia.cingolani@uniba.it}{\textcolor{black}{silvia.cingolani@uniba.it}}}
		\\ 
		\\Marco Gallo 
		\orcidicon{0000-0002-3141-9598}
		\\ \normalsize{Dipartimento di Matematica e Fisica}
		\\ \normalsize{Universit\`{a} Cattolica del Sacro Cuore}
		\\ \normalsize{Via della Garzetta 48, 25133 Brescia, Italy}
		\\ \normalsize{\href{mailto:marco.gallo1@unicatt.it}{\textcolor{black}{marco.gallo1@unicatt.it}}}
		\\ 
		\\Kazunaga Tanaka 
		\orcidicon{0000-0002-1144-1536}
		\\ \normalsize{Department of Mathematics}
		\\ \normalsize{School of Science and Engineering}
		\\ \normalsize{Waseda University}
		\\ \normalsize{3-4-1 Ohkubo, Shijuku-ku, Tokyo 169-8555, Japan}
		\\ \normalsize{\href{mailto:kazunaga@waseda.jp}{\textcolor{black}{kazunaga@waseda.jp}}}
	}

\date{}

\maketitle

%
%
%

\begin{abstract} 
In this paper we study the following nonlinear fractional Choquard-Pekar equation 
\begin{equation}\label{eq_abstract}
(-\Delta)^s u + \mu u =(I_\alpha*F(u)) F'(u) \quad \hbox{in}\ \mathbb{R}^N, \tag{$*$}
\end{equation}
where $\mu>0$, $s \in (0,1)$, $N \geq 2$, $\alpha \in (0,N)$, $I_\alpha \sim \frac{1}{|x|^{N-\alpha}}$ is the Riesz potential, and $F$ is a general subcritical nonlinearity. The goal is to prove existence of multiple (radially symmetric) solutions $u \in H^s(\mathbb{R}^N)$, by assuming $F$ odd or even: we consider both the case $\mu>0$ fixed and the case $\int_{\mathbb{R}^N} u^2 =m>0$ prescribed.
Here we also simplify some arguments developed for $s=1$ \cite{CGT4}.

 A key point in the proof is given by the research of suitable multidimensional odd paths, which was done in the local case by Berestycki and Lions \cite{BL2}; for \eqref{eq_abstract} the nonlocalities play indeed a special role.
In particular, some properties of these paths are needed in the asymptotic study (as $\mu$ varies) of the mountain pass values of the unconstrained problem, then exploited to describe the geometry of the constrained problem and detect infinitely many normalized solutions for any $m>0$.

The found solutions satisfy in addition a Pohozaev identity: in this paper we further investigate the validity of this identity for solutions of doubly nonlocal equations under a $C^1$-regularity.
\end{abstract}


\medskip




{
  \hypersetup{linkcolor=black} 
  \tableofcontents
}



\section{Introduction
}
\label{sec_intro_choquard}

Given a nonlinearity $F \in C^1(\R,\R)$ and set $f:=F'$, we are interested to seek for multiple solutions of the nonlocal system 
\begin{equation}\label{eq_SN_syst}
\parag
{(-\Delta)^s u + \mu u = \phi f(u) \quad \hbox{ in $\R^N$}, \\
(-\Delta)^{\alpha/2} \phi = F(u) \quad \hbox{ in $\R^N$},
}
\end{equation}

where $N\geq 2$, $s \in (0,1)$ and $\alpha\in (0,N)$. 
In literature this is system has several physical motivations and it is usually called \emph{fractional Schr\"odinger-Newton system}.

By solving the second equation of \eqref{eq_SN_syst} through the Riesz-Poisson kernel $I_{\alpha}$, and substituting $\phi$ in the first, we come up with the single equation
\begin{equation}\label{eq_Choquard_genericaF}
(- \Delta)^s u + \mu u =(I_\alpha*F(u)) f(u) \quad \hbox{in}\ \R^N, 
\end{equation}
which is also known as \emph{fractional} \emph{Hartree} (or \emph{Choquard-Pekar}) equation.


\smallskip

Applications of \eqref{eq_Choquard_genericaF} can be found in quantum chemistry \cite{ArMe,DOS,HMT0} (see also \cite{CFHMT} for some orbital stability results): here \eqref{eq_Choquard_genericaF} appears in the study of the mean field limit of weakly interacting multi-particle systems.
In particular the equation applies to the study of graphene \cite{LuMoMu}, where the nonlocal nonlinearity describes the short time interactions between particles. 
Doubly nonlocal equations appear also in the dynamics of populations \cite{CDV}. 

One of the main applications anyway arises in the study of exotic stars: minimization properties related to \eqref{eq_Choquard_genericaF} play indeed a fundamental role in the mathematical description of the dynamics of pseudorelativistic boson stars and their gravitational collapse, as well as the evolution of attractive fermionic systems, such as white dwarf stars \cite{ElSc0, FroL1, FTY, FJL1, HS0, Len1, LT0, LiYa0, MiSc0, Ngu1}.
In fact, the nonlocality of the operator is expression of the relativistic velocity of the particles, while the nonlocal source describes the mutual interaction of the particles in the same body. Here the study of the ground states to \eqref{eq_Choquard_genericaF} gives information on the size of the critical initial conditions for the solutions of the corresponding pseudorelativistic equation, where a critical value is given by the Chandrasekhar \cite{Chd0} limiting mass. 
In particular, when $s=\frac{1}{2}$, $N=3$, $\alpha=2$ and $f(u)=u$, we obtain \cite{FraL,LeLe,HL0}
\begin{equation}\label{eq_boson_massless}
\sqrt{-\Delta}u + \mu u = \left(\frac{1}{4 \pi |x|}*u^2\right) u \quad \hbox{in $\R^3$}.
\end{equation}
We recall that equation \eqref{eq_boson_massless}, in the case of non-relativistic Laplacian $-\Delta$, was elaborated in 1954 by Pekar \cite{Pek0} (see also \cite{LR0}), and subsequently considered as a physical model by many authors \cite{Cho0, Lie1, Har0, Pen1, TM0}; the fractional Laplacian operator was instead introduced by Feynman \cite{Fey0} in 1948 (see also \cite{Las1}).

\begin{Remark}
In \eqref{eq_boson_massless} $f(\sigma )=|\sigma|^{r-2}\sigma$ with $r=2$ is $L^2$-critical: in this paper, when dealing with the mass-constrained problem, we essentially address the subcritical case $r\in (\frac{5}{3}, 2)$ (see condition \hyperref[(CF3)]{\textnormal{(CF3)}} below), but we believe that this result, together with the developed minimax tools, can be a first step towards the study of the $L^2$-mass critical (and supercritical) case, since for these problems the minimization approach is generally not well posed. 
\end{Remark}

If $u$ is a solution of \eqref{eq_Choquard_genericaF}, then the wave function $\psi(x,t) = e^{i \mu t} u(x)$
is a solitary wave of the time-dependent equation (when, for example, $f(\psi)=g(|\psi|)$ or $f(\psi)=g(|\psi|)\psi$, and similarly $F$) 
\begin{equation} \label{eq:1.3}
i \psi_t = (- \Delta)^s \psi - \big( I_{\alpha}*F(\psi)\big) f(\psi) \quad \text{in $\R^N \times (0, +\infty) $}.
\end{equation}
The study of standing waves of \eqref{eq:1.3} has been pursed in two main directions, which opened two different challenging research fields.

A first topic regards the search for solutions of \eqref{eq_Choquard_genericaF} with a prescribed frequency $\mu$ and free mass, the so-called \emph{unconstrained} problem. 
The second line of investigation of the problem \eqref{eq_Choquard_genericaF} consists of prescribing the mass $m >0$ of $u$, thus conserved by $\psi$ in time
$$
\int_{\R^3} |\psi(x,t)|^2 \, dx= m \quad \forall \, t \in [0, +\infty),
$$
and letting the frequency $\mu$ to be free. Such problem is usually said \emph{constrained} problem, and its study is also strictly related to the dynamical properties of the solutions of \eqref{eq:1.3} \cite{Geo0}.

\medskip

For the unconstrained problem, the first investigations on the existence of solutions when $s=1$ go back to \cite{Lie2,Men0,Stu1,CSV,MPT,MZo0};
variational methods were also employed to derive existence and qualitative results of standing wave solutions for more generic values of $\alpha \in (0,N)$ and of pure power nonlinearities $F(\sigma )= \frac{1}{p} |\sigma|^p$ by Moroz and Van Schaftingen \cite{MS0} (see also \cite{MS3}), and then generalized for doubly nonlocal equations in the papers \cite{DSS1} by D'Avenia, Siciliano and Squassina. 
In particular the authors in \cite{DSS1} considered the special model
\begin{equation}\label{eq:1.5}
(- \Delta)^s u + \mu u =(I_\alpha \ast |u|^p) |u|^{p-2} u \quad \hbox{in}\ \R^N, 
\end{equation}
and they proved that \eqref{eq:1.5} has solutions if 
\begin{equation}\label{eq:1.6}
2^{\#}_{\alpha}:=\frac{N+\alpha}{N} < p < \frac{N+\alpha}{N-2s} =: 2^*_{\alpha,s}.
\end{equation}
When dealing with variational (and regular) solutions, they proved that range \eqref{eq:1.6} is optimal. 
Other results can be found in \cite{Luo0} for general superlinear nonlinearities, in \cite{ChLi} for the non-autonomous case, in \cite{LMZ} for critical case, in \cite{LLL} for zero mass case.

Recently in \cite{CGT3} the authors considered the problem \eqref{eq_Choquard_genericaF} when $F$ is a Berestycki-Lions \cite{BL1} type function under the following general assumptions, extending to the fractional case what has been done by Moroz and Van Schaftingen \cite{MS2} 
 \begin{itemize}
\item[(F1)] \label{(F1)}
$F \in C^1(\R, \R)$;
\item[(F2)] \label{(F2)}
there exists $C >0$ such that, for every $\sigma \in \R$, $$|\sigma f(\sigma )| \leq C \big(|\sigma|^{2^{\#}_{\alpha}} + |\sigma|^{2^*_{\alpha,s}}\big);$$
\item[(F3)] \label{(F3)}
$$\lim_{\sigma \to 0} \frac{F(\sigma )}{|\sigma|^{2^{\#}_{\alpha}}} =0, \quad
\lim_{ \sigma \to + \infty} \frac{F(\sigma )}{|\sigma|^{2^*_{\alpha,s}}} =0;$$
\item[(F4)] \label{(F4)}
$F\not\equiv 0$, that is, there exists $\sigma_0 \in \R$, $\sigma_0 \neq 0$ 
 such that $F(\sigma _0) \neq0$.
 \end{itemize}

In particular they prove the existence of a Pohozaev minimum, which is equivalently proved to have a mountain pass structure, satisfying the Pohozev identity \cite[Theorem 1.1]{CGT3}. Further qualitative properties, such as regularity, summability, positivity and asymptotic decay, have been then investigated in \cite{CGT3, CG1,Gal1} (see also \cite{GalT}).

The assumptions \hyperref[(F1)]{\textnormal{(F1)}}--\hyperref[(F4)]{\textnormal{(F4)}}, which rely essentially only on the growth of the nonlinearity in zero and at infinity -- not including regularity, homogeneity, Ambrosetti-Rabinowitz-type or monotonicity conditions -- can be considered, from a variational point of view, \emph{almost optimal}. 
They include for instance the most common power type functions $f(\sigma) \sim |\sigma|^{p-1} \sigma$, $f(\sigma) \sim |\sigma|^{p}$, but also combined powers representing cooperation $f(\sigma) \sim |\sigma|^{p-1} \sigma +|\sigma|^{q-1}\sigma$, $f(\sigma) \sim |\sigma|^p+|\sigma|^q$ and competition $f(\sigma) \sim |\sigma|^{p-1} \sigma-|\sigma|^{q-1}\sigma $, $f(\sigma) \sim |\sigma|^p-|\sigma|^q$, as well as asymptotically linear 
sources $f(\sigma) \sim \frac{\sigma^3}{1+\sigma^2}$, $f(\sigma) \sim \frac{|\sigma|^3}{1+\sigma^2}$ typical of saturation effects arising in nonlinear optics \cite{DLQWZPH, 
RAAY}, 
and many others.

\smallskip

The existence of an infinite number of standing wave solutions to \eqref{eq:1.5} was faced by \cite{DSS1} (while previous results for $s=1$ go back to \cite{Lio1,ClSa}): here the homogeneity of the source plays a crucial role, and similar ideas have been applied for $s=1$ in \cite{Ack0,QRT} in presence of more general sources satisfying Ambrosetti-Rabinowitz type conditions. 
We highlight that all the abovementioned papers work with odd nonlinearities $f$. 

It remains open the existence of infinitely many radially symmetric solutions for the nonlinear fractional Choquard equation \eqref{eq_Choquard_genericaF} under the optimal assumptions \hyperref[(F1)]{\textnormal{(F1)}}--\hyperref[(F4)]{\textnormal{(F4)}} and symmetric conditions on the nonlocal source term $(I_\alpha*F(u))f(u)$, that is
\begin{itemize}
\item[(F5)] \label{(F5)}
$F$ is odd or even.
\end{itemize}

Differently from existence \cite{CGT3}, the symmetry of $F$ is crucial for the multiplicity of solutions.
Here we assume $F$ to be \emph{odd} or \emph{even}, which guarantees the evenness of the energy functional associated to \eqref{eq_Choquard_genericaF}. 
Contrary to the source-local case \cite{BL2,HT0,Iko2,CGT1}, where the nonlinear term is usually assumed odd, for nonlocal sources we emphasize that both the cases $F$ odd and even are meaningful (see also \cite{DMP} where \hyperref[(F5)]{\textnormal{(F5)}} is assumed to achieve existence of a single nonradial solution).

By virtue of \cite{Pal0}, radially symmetric solutions to \eqref{eq:1.7} can be searched as critical points of the $C^1$-functional $\mc{J}_{\mu}: H^s_r(\R^N) \to \R$ 
$$
\mc{J}_{\mu} (u):=\half \int_{\R^N} |(-\Delta)^{s/2} u|^2 \, dx + \frac{\mu}{2} \int_{\R^N} u^2 \, dx 
-\half \int_{\R^N} (I_\alpha*F(u))F(u)\, dx,
$$
where $H^s_r(\R^N)$ denotes the fractional Sobolev space of radially symmetric functions. 
We establish thus the following result.
\begin{Theorem}\label{S:1.3}	
Suppose $N\geq 2$, $\alpha \in (0, N)$, $s\in (0,1)$ and $\mu>0$ be fixed. Assume that \hyperref[(F1)]{\textnormal{(F1)}}--\hyperref[(F5)]{\textnormal{(F5)}} hold. Then there exist countably many radial solutions $(u_n)_{n}$ of \eqref{eq_Choquard_genericaF} which satisfy the Pohozaev identity
 \begin{equation}\label{eq:2.5}
\frac{N-2s}{2} \int_{\R^N}|(-\Delta)^{s/2} u|^2 + \frac{N }{ 2} \mu \int_{\R^N} u^2 - \frac{N + \alpha}{2} \int_{\R^N} \big(I_{\alpha}*F(u)\big) F(u)=0.
\end{equation}
 Moreover we have 
$$\mathcal{J}_{\mu}(u_n) \to +\infty \quad \hbox{as $ n \to + \infty$}.$$
\end{Theorem}
Our multiplicity result is the counterpart of what done in \cite{BL2} (for $N\geq 3$, see also \cite{BGK} for $N=2$) in the local case with odd nonlinearities and extend the existence result in \cite{DSS1}. 

Due to the generality of the source, in order to get compactness, we use a variant of the Palais-Smale condition \cite{HT0,IT0}, which takes into account the Pohozaev identity, through which we prove a suitable deformation theorem in an augmented space \cite{CT2}, and eventually we use the genus tool to get infinitely many minimax solutions.

In order to implement minimax arguments, anyway, a key point is the construction of multidimensional odd paths. 
When $f$ satisfies some Ambrosetti-Rabinowitz condition (i.e., $F$ can be estimated from below by an homogeneous function $|\sigma|^p$), the construction of such a path classically relies on the equivalence of the $H^s$-norm and the $L^p$-norm on finite dimensional subspaces of $H^s(\R^N)$. When such condition is no longer available, a finer construction is needed: in the celebrated paper \cite{BL2} Berestycki and Lions build this path for a local problem by exploiting an inductive process based on piecewise affine functions.

In our nonlocal case, in order to prove the existence of multiple solutions, unlike the elaborated approach of \cite{BL2} we can obtain the existence of a multidimensional odd path by exploiting the positivity of the Riesz potential functional. 
This gives also a simplified approach, differently from what was presented in \cite{CGT4}. See anyway Remark \ref{rem_general_kernels} below for some comments.

\medskip

We move now to study the constrained case, which appears more delicate. 
The existence of $L^2$-normalized solutions was investigated when $F(\sigma)=\frac{1}{p}|\sigma|^p$ in \cite{Wu0,GS1} (see also \cite{CFHMT} where some Ambrosetti-Rabinowitz conditions are assumed and \cite{GZ0} where scattering properties are investigated for $L^2$-supercritical evolutive problems).
More recently, the authors in \cite{CGT2} obtained existence of a solution $u\in H^s_r(\R^N)$ to
\begin{equation}\label{eq:1.7}
\left \{
\begin{aligned}
& (- \Delta)^s u + \mu u =(I_\alpha*F(u))f(u) \quad \hbox{in}\ \R^N, \\ 
&\int_{\R^N} u^2 dx = m, \; \; \mu>0, 
\end{aligned}
\right. 
\end{equation}
assuming that $F$ satisfies \hyperref[(F1)]{\textnormal{(F1)}}, \hyperref[(F4)]{\textnormal{(F4)}} and it is $L^2$-subcritical, namely
\begin{itemize}
\item[(CF2)] \label{(CF2)}
there exists $C >0$ such that, for every $\sigma \in \R$, $$|\sigma f(\sigma )| \leq C \big(|\sigma|^{2^{\#}_{\alpha}} + |\sigma|^{2^m_{\alpha,s}}\big);$$
\item[(CF3)] \label{(CF3)}
$$\lim_{\sigma \to 0} \frac{F(\sigma )}{|\sigma|^{2^{\#}_{\alpha}}} =0, \quad
\lim_{ \sigma \to + \infty} \frac{F(\sigma )}{|\sigma|^{2^m_{\alpha,s}}} =0,$$
\end{itemize}
where
$$2^m_{\alpha,s}:= \frac{N+ \alpha+2s}{N}.$$

When $s=1$, multiplicity of radial standing wave solutions to \eqref{eq:1.7} with prescribed $L^2$-norm has been faced for powers again by Lions in \cite{Lio1} (see also \cite{CJ0} for the planar logarithmic Choquard equation). 
As regards instead the case of general nonlinearities $f$, recently Bartsch et al. \cite{BaLiLi} obtained the existence of infinitely many solutions of \eqref{eq:1.7} by assuming that $f$ is an odd function which satisfies monotonicity and Ambrosetti-Rabinowitz conditions; 
here the authors in \cite{BaLiLi} rely on concentration-compactness arguments, together with the use of a stretched functional.
When $s\in (0,1)$ some results have been achieved in \cite{YZZ} in presence of local perturbations, and in \cite{LiLuo} for the supercritical case.
In this paper we address the case of general $f$, when monotonicity and Ambrosetti-Rabinowitz type conditions do not hold, or $f$ is not odd.

A typical approach is to characterize solutions to \eqref{eq:1.7} as critical points of the $C^1$-functional $\mc{L}: H^s_r(\R^N) \to \R$
$$ \mathcal{L}(u) := \half \int_{\R^N} |(-\Delta)^{s/2} u|^2\, dx -\half \int_{\R^N} (I_\alpha*F(u))F(u)\, dx,$$
constrained on the sphere 
$$ \mathcal{S}_m := \left \{ u \in H^s_r(\R^N) \mid \int_{\R^N} u^2 \, dx= m \right\},$$
and find here a minimum of this functional (whenever bounded). 
On the other hand, this approach requires much attention in order to control the tails of the functions. That is why, in the spirit of \cite{HT0}, we work instead with a Lagrangian formulation of the nonlocal problem \eqref{eq:1.7}.
We highlight the advantage of this method, that can be suitably adapted to derive multiplicity results of normalized solutions in several different frameworks, and does not require boundedness of the functional on the $L^2$-ball.

Through this formulation we prove the following result.
\begin{Theorem}\label{S:1.1_choq}	
Suppose $N\geq 2$, $\alpha \in (0, N)$, $s\in (0,1)$ and \hyperref[(F1)]{\textnormal{(F1)}}-\hyperref[(CF2)]{\textnormal{(CF2)}}-\hyperref[(CF3)]{\textnormal{(CF3)}}-\hyperref[(F4)]{\textnormal{(F4)}}-\hyperref[(F5)]{\textnormal{(F5)}}.
\begin{itemize}
\item[(i)] For any $k \in \N$ there exists $m_k \geq 0$ such that for every $m > m_k$, the problem \eqref{eq:1.7} has at least $k$ pairs of nontrivial, distinct, radially symmetric solutions $(\mu_n, \pm u_n)_{n=1\dots k}$, which satisfy the Pohozaev identity \eqref{eq:2.5}. 
\item[(ii)] Assume in addition an $L^2$-subcritical growth also at zero, i.e.
\begin{itemize}
\item[\textnormal{(CF4)}] \label{(CF4)} 
$$	
\lim_{\sigma \to 0} \frac{|F(\sigma )|}{|\sigma|^{2^m_{\alpha,s}}} = + \infty;
$$
additionally, if $F$ is odd, assume that there exists $\delta_0>0$ such that $F$ has a constant sign in $(0,\delta_0]$ and
\begin{equation}\label{eq_cond_more_gen}
\sup_{\sigma \in (0, \delta_0],\, h\in [0,1]} \frac{F(\sigma h)}{F(\sigma )}< +\infty;
\end{equation}
for example, this is satisfied if $\abs{F}$ is assumed non-decreasing in $[0,\delta_0]$. 
\end{itemize}
Then $m_k=0$ for each $k \in \N$, that is for any $m >0$ the problem \eqref{eq:1.7} has countably many pairs of solutions $(\mu_n, \pm u_n)_{n}$ satisfying the Pohozaev identity \eqref{eq:2.5}. Moreover $\mathcal{L}(u_n) <0$, $n \in \N$ and we have 
$$\mathcal{L}(u_n) \to 0 \quad \hbox{as $ n \to + \infty$}.$$
\end{itemize}
\end{Theorem}

\begin{Remark}
\label{rem_extra_cond}
We comment condition \eqref{eq_cond_more_gen}. Set 
$$ M:=\sup_{\sigma \in (0, \delta_0],\, h\in [0,1]} \frac{F(\sigma h)}{F(\sigma )}< +\infty$$
we have, when $\abs{F}$ is non-decreasing, $M=1$. As a nontrivial example one can consider $\beta \in (2^{\#}_{\alpha},2^m_{\alpha,s})$ and $F$ oscillating near zero between $|\sigma|^{\beta}$ and $2|\sigma|^{\beta}$, so that $M\leq 2$; for instance the odd extension of
$$F(\sigma ):= \sigma^{\beta} \big(2+ \sin(\tfrac{1}{\sigma})\big) \quad \hbox {as $\sigma \to 0^+$}.$$
If instead $F$ oscillates (not strictly) between $|\sigma|^{\beta_1}$ and $|\sigma|^{\beta_2}$, with $2^{\#}_{\alpha} < \beta_1 < \beta_2 < 2^m_{\alpha,s}$, then $M=+\infty$; thus for instance the odd extension of
$$F(\sigma ):= \sigma^{\beta_1} \big(1+ \sin(\tfrac{1}{\sigma})\big) +\sigma^{\beta_2} \big(1- \sin(\tfrac{1}{\sigma})\big) \quad \hbox {as $\sigma \to 0^+$}$$
is not covered by \eqref{eq_cond_more_gen}.
\end{Remark}

As mentioned before, in the spirit of \cite{CGT1, CGT2, CGT4}, the approach to the proof of this result relies on some Lagrangian formulation, developed on a Pohozaev mountain geometry and a Palais-Smale-Pohozaev compactness condition: this setting allows to detect a minimax structure on a product space and thus the existence of infinitely many solutions.

Of key importance here is the asymptotic analysis of the mountain pass values of the unconstrained problem (as $\mu\to 0^+$ and $\mu \to +\infty$). 
To this aim, again, a delicate issue is the construction of suitable multidimensional odd paths. 
To prove the existence of multiple solutions for $m\gg 0$ (point ($i)$ of Theorem \ref{S:1.1_choq}) we can apply the same arguments of Theorem \ref{S:1.3}. 
A similar approach can be implemented also to gain existence of infinitely many solutions for any $m>0$ when $F$ is even (first part of point $(ii)$ of Theorem \ref{S:1.1_choq}), since in this case $F$ can be assumed nonnegative in a neighborhood of the origin. 
See anyway Remark \ref{rem_general_kernels} below.

A quite delicate issue, instead, comes up when $F$ is odd and $m>0$ is arbitrary. 
As in \cite{CGT4}, differently from \cite{HT0,CGT1}, we need to implement a new approach to gain the existence of an admissible odd path which satisfies some finer condition. 
Such a path is built by using suitable annuli: the interactions, produced by the Riesz potential, between the characteristic functions of these annuli have to be specifically investigated, and some sharp estimates for the Riesz potential \cite{Thim0} are here used to this goal. 
The power $\alpha \in (0,N)$ describes the strength of this interaction, and it reveals to be stronger when $\alpha$ is small, i.e. $\alpha\in (0,1]$.

\begin{Remark}\label{rem_general_kernels}
We highlight that the easier approach for building a multidimensional path, based on the positivity of the Riesz kernel (Proposition \ref{prop_positivity_Riesz}), cannot generally be applied to more generally frameworks (also if $F$ is even); for examples, when dealing with kernels $K=K(x,y)$ which does not make the functional
$$g \mapsto \int_{\R^N} \int_{\R^N} K(x,y) g(x) g(y) dx dy$$
 positive (e.g., $K(x,y)$ sign-changing): some examples are given by $K(x,y)=\log(|x-y|)$, $K(x,y)= (\pm I_{\alpha}\chi_{B_1} \mp I_{\beta} \chi_{B_1^c})(x-y)$ for some $\alpha,\beta \in (0,N)$, but also $K(x,y)= \frac{||x|-|y|||}{|x||y|}$. 
In this case, the approach here developed, based on suitable annuli, might instead be adapted. 
This is an interesting line of research for the future. 
\end{Remark}

Finally, in this paper we deal with the Pohozaev identity, which essentially says that 
$$\frac{d}{d\theta} \mc{J}_{\mu}(u(\cdot/e^{\theta}))_{|\theta=0}=0.$$
This identity may be useful in different occasions: nonexistence results, information on multipliers, semiclassical analysis (see \cite{CG0}), and so on. 
We note that in Theorems \ref{S:1.3} and \ref{S:1.1_choq} we find solutions satisfying the Pohozaev identity; however, it is not known whether this identity holds for all weak solutions or not. 
When $f$ is a power, such identity has been proved \cite[equation (6.1)]{DSS1}: in this case the solutions are indeed $C^2$, which is a key point in the proof.
In this paper, in the spirit of \cite{BKS} and \cite{Dje0}, we show this identity by asking only a $C^1$-regularity; we highlight that in the proof we avoid the use of the $s$-harmonic extension. 
We get thus the following result.

\begin{Theorem}
\label{thm_pohoz_frac_choq}
Let $u\in H^s(\R^N)\cap L^{\infty}(\R^N)$ be a weak solution of \eqref{eq_Choquard_genericaF}, and assume \hyperref[(F1)]{\textnormal{(F1)}}-\hyperref[(F2)]{\textnormal{(F2)}}.
Assume moreover $f\in C^{0,\sigma}_{loc}(\R)$ for some $\sigma \in (0,1]$ and one of the following:
\begin{itemize}
\item $s \in [\frac{1}{2}, 1)$, 
\item $s \in [\frac{1}{4}, \frac{1}{2})$, $\alpha \in ( 1-2s, N)$ and $\sigma \in (\frac{1-2s}{2s}, 1]$,
\item $s \in (0, \frac{1}{2})$, $\alpha \in (0,2)$ and $\sigma \in (1-2s, 1]$.
\end{itemize}
Then $u \in C^{1,\gamma}(\R^N)$ for some $\gamma \in (\max\{0,2s-1\},1) $, and $u$ satisfies the Pohozaev identity \eqref{eq:2.5}, or equivalently
$$\frac{1}{2^*_{\alpha,s}} \|(-\Delta)^{s/2} u\|_2^2 + \frac{\mu}{2^{\#}_{\alpha}} \|u\|_2^2 - \mc{D}(u)=0$$
where $\mc{D}(u) := \int_{\R^N} (I_\alpha*F(u))F(u)$. The result in particular applies to nonnegative weak solutions $u \in H^s(\R^N)$ of \eqref{eq_Choquard_genericaF}.
\end{Theorem}

\smallskip

The paper is organized as follows. 
In Section \ref{sec_prelim} we recall some facts about the fractional Laplacian and the Riesz potential.
Section \ref{section:4..1} is dedicated to the construction of different multidimensional paths, both with a simple approach and a more refined approach, suitable for existence of infinitely many normalized solutions.
Then these refined paths are exploited in some asymptotic analysis in Section \ref{sec_choquard_asympt_MP}, and these results are then applied to achieve multiplicity results in Section \ref{sec_exist_sol}.
Finally, in Section \ref{sec_Pohozaev}, we prove a Pohozaev identity for the fractional Choquard equation.


\section{Preliminaries}
\label{sec_prelim}

Set, for $u\in L^2(\R^N)$ with $|\xi|^{2s} \mc{F}(u) \in L^2(\R^N)$, the fractional Laplacian \cite{DnPV}
$$
(-\Delta)^s u := \mc{F}^{-1}(|\xi|^{2s} \mc{F}(u));
$$
when $u$ is regular enough, we can write
$$(-\Delta)^s u(x) = C_{N,s} \textnormal{PV} \int_{\R^N} \frac{u(x)-u(y)}{|x-y|^{N+2s}} dy, \quad x \in \R^N$$
where $C_{N,s}:=\frac{4^s \Gamma(\frac{N+2s}{2})}{\pi^{N/2} |\Gamma(-s)|}>0$. 
We set the homogeneous fractional Sobolev space
$$ D^{s,2}(\R^N):= \left\{ u \hbox{ measurable} \mid (-\Delta)^{s/2} u \in L^2(\R^N) \right \}$$
and the fractional Sobolev space
$$H^s(\R^N) := L^2(\R^N) \cap D^{s,2}(\R^N) $$ 
endowed with
$\norm{u}_{H^s}^2 := \norm{u}_2^2 + \norm{(-\Delta)^{s/2}u}_2^2.$
An equivalent form of the seminorm is given by
$$
[u]_{H^s(\R^N)}^2 := \frac{C_{N,s}}{2} \int_{\R^N} \int_{\R^N} \frac{|u(x)-u(y)|^2}{|x-y|^{N+2s}} dx dy 
= \norm{(-\Delta)^{s/2} u}_2^2.
$$
Consider moreover the subspace of radially symmetric functions
$$H^s_r(\R^N) := \big\{ u \in H^s(\R^N) \mid u(x)= u(|x|) \big\}$$
and recall that $H^s(\R^N) \hookrightarrow L^p(\R^N)$ 
for every $p \in [2, 2^*_s]$, and $H^s_r(\R^N) \hookrightarrow \hookrightarrow L^p(\R^N)$ for every $p \in (2, 2^*_s)$.

A sufficient condition in order to have $(-\Delta)^s u$ well defined pointwise is given by \cite[Proposition 2.4]{Silv0} (see also \cite[Proposition 2.15]{Gar0} and \cite[Lemma 2.4]{BKS}). We use the notation
$$C^{\gamma}(\R^N):= C^{[\gamma], \gamma-[\gamma]}(\R^N)$$
and $Lip(\R^N):=C^{0,1}(\R^N)$, and similarly $C^{\gamma}_{loc}(\R^N)$, $Lip_{loc}(\R^N)$, $C^{\gamma}_{c}(\R^N)$ and $Lip_{c}(\R^N)$.
\begin{Proposition}
\label{prop_well_posed}
Let $s\in (0,1)$. Assume 
\begin{equation} \label{eq_spazio_gener_s}
\int_{\R^N} \frac{|u(x)|}{(1+|x|)^{N+2s}} dx< \infty;
\end{equation}
e.g., $u \in (L^1+ L^{\infty})(\R^N)$. 
Assume moreover $u \in C^{\gamma}_{loc}(\R^N)$ for some $\gamma >2s$.
Then
\begin{align*}
(-\Delta)^s u(x) &= C_{N,s} \lim_{\eps \to 0} \int_{\R^N \setminus B_{\eps}(x)} \frac{u(x)-u(y)}{|x-y|^{N+2s}} dy \\
&= \frac{C_{N,s}}{2} \int_{\R^N}\frac{ 2u(x)-u(x+y)-u(x-y)}{|y|^{N+2s}} dy
\end{align*}
for every $x \in \R^N$, and in this case we have $(-\Delta)^s u \in C(\R^N)$. Moreover the last integral is absolutely convergent.
\end{Proposition}

\claim Proof.
We check only the absolute convergence. Indeed, let $x \in \R^N$ and $R>2|x|+1$. Notice that, for $|y|\geq R$, we have, for $|y|\geq R$,
$$|x+y| \geq |y|-|x| \geq R - |x|>|x|+1 $$
and
$$|x+y| - |x| \geq \frac{|x+y|+1}{|x|+2} $$
thus
\begin{align*}
\int_{B_R^c}\frac{ |2u(x)-u(x+y)-u(x-y)|}{|y|^{N+2s}} dy &\leq \int_{B_R^c}\frac{2 |u(x)|}{|y|^{N+2s}} dy + \int_{B_R^c}\frac{|u(x+y)|}{|y|^{N+2s}}dy + \int_{B_R^c}\frac{|u(x-y)|}{|y|^{N+2s}} dy \\
&\leq 2 |u(x)| \int_{B_R^c}\frac{1}{|y|^{N+2s}} dy + 2\int_{B_{|x|+1}^c}\frac{|u(z)|}{(|z|-|x|)^{N+2s}}dz\\ 
&\leq C_R |u(x)| + 2(2+|x|)^{N+2s}\int_{B_{|x|+1}^c}\frac{|u(z)|}{(|z|+1)^{N+2s}}dz < \infty.
\end{align*}
Let now $s \in (0, \frac{1}{2})$. Then, being $u \in C^{0,\gamma}_{loc}(\R^N)$ for some $\gamma>2s$,
$$\int_{B_R}\frac{ |2u(x)-u(x+y)-u(x-y)|}{|y|^{N+2s}} dy \leq 2 C \int_{B_R}\frac{ 1}{|y|^{N+2s-\gamma}} dy < \infty;$$
notice that a similarly argument shows also that the first integral of the Proposition's claim does not need the principal value, being absolute convergent.

If instead $s \in [\frac{1}{2}, 1)$, then, being $u \in C^{1,\gamma}_{loc}(\R^N)$ for some $\gamma>2s-1$, for each $x,y \in \R^N$ there exists $\sigma=\sigma(x,y) \in (0,1)$ such that
\begin{align*}
\int_{B_R}\frac{ |2u(x)-u(x+y)-u(x-y)|}{|y|^{N+2s}} dy & = \int_{B_R}\frac{ |\nabla u(x+\sigma y)\cdot y - \nabla u(x-\sigma y) \cdot y|}{|y|^{N+2s}} dy \\
& \leq \int_{B_R}\frac{2\sigma}{|y|^{N+2s-\gamma-1}} dy < \infty.
\end{align*}
Joining the pieces, we have the claim.
\QED

\bigskip

Set then the Riesz potential
$$
I_{\alpha}(x) := \frac{C_{N,\alpha}}{|x|^{N-\alpha}}
$$
where $C_{N,\alpha}:=\frac{\Gamma(\frac{N-\alpha}{2})}{2^{\alpha} \pi^{N/2} \Gamma(\frac{\alpha}{2})}>0$, and define
$$\mc{D}_{\alpha}(g,h):= \int_{\R^N} (I_{\alpha}*g)h = \int_{\R^N} \int_{\R^N} \frac{g(x) h(y)}{|x-y|^{N-\alpha}} dx dy.$$
The following theorem ensures the well posedness of the Riesz potential (see \cite[Theorem 4.3]{LiLo} and \cite[pages 61-62]{Lan0}). 
\begin{Proposition}[Hardy-Littlewood-Sobolev inequality]\label{prop_HLS}
Let $\alpha \in (0,N)$.
\begin{itemize}
\item Let $g$ be a measurable function. Then $I_{\alpha}*g$ is finite almost everywhere if and only if
\begin{equation}\label{eq_buona_def_riesz}
\int_{\R^N} \frac{|g(x)|}{(1+|x|)^{N-\alpha}}dx <\infty.
\end{equation}
In particular, $I_{\alpha}*g$ is well defined if $g \in L^1_{loc}(\R^N) \cap L^r(B_R^c)$ for some $R\geq 0$ and some $r \in [1, \frac{N}{\alpha})$.
Moreover, if \eqref{eq_buona_def_riesz} does not hold, then $I_{\alpha}*|g| \equiv \infty$. 
\item 
Let $r \in [1, \frac{N}{\alpha})$. Then the map
$$g \in L^r(\R^N) \mapsto I_{\alpha}*g \in L^{\frac{Nr}{N-\alpha r}}(\R^N)$$
is continuous.
 In particular, since the operator is linear,
$$g_n \wto g \hbox{ weakly in $L^r(\R^N)$} \implies I_{\alpha}*g_n \wto I_{\alpha}*g \hbox{ weakly in $L^{\frac{Nr}{N-\alpha r}}(\R^N)$}.$$
\item Let $r,t \in (1,+\infty)$ be such that $\frac{1}{r}+\frac{1}{t}= \frac{N+\alpha}{N}$. Then there exists a constant $C=C(N,\alpha,r,t)>0$ such that
$$\pabs{\mc{D}_{\alpha}(g,h)}\leq C \norm{g}_r \norm{h}_t$$
for all $g \in L^r(\R^N)$ and $h \in L^t(\R^N)$.
\end{itemize}
\end{Proposition}

We observe the following: if $g \in \mc{S}$ \cite[Lemma 5.1.2]{Ste0} or if $\alpha \in (0, \frac{N}{2})$ and $g \in L^{\frac{2N}{N+2\alpha}}(\R^N)$ \cite[Corollary 5.10]{LiLo} then we have 
$$\mc{D}_{\alpha}(g,g)=\int_{\R^N} \big(I_{\alpha}*g\big) g dx = \int_{\R^N} \mc{F}\big(I_{\alpha}*g\big) \mc{F}(g) d\xi= \int_{\R^N} \mc{F}(I_{\alpha}) |\mc{F}(g)|^2 d\xi = \int_{\R^N} \frac{|\mc{F}(g)|^2}{|\xi|^{\alpha}} d\xi \geq 0$$
This shows that
$$g \mapsto \mc{D}_{\alpha}(g,g)$$ 
is a positive bilinear form. 
A more general result can be adapted from $\alpha=2$ \cite[Theorem 9.8]{LiLo} to a generic $\alpha \in (0,N)$ as follows.
\begin{Proposition}[\cite{LiLo}]\label{prop_positivity_Riesz}
Let $g:\R^N \to \R$ measurable be such that
$$\mc{D}_{\alpha}(|g|,|g|)< \infty.$$
Then
$$\mc{D}_{\alpha}(g,g)\geq 0 $$
and the above quantity is zero if and only if $g\equiv 0$ almost everywhere.
In particular the following representation holds 
$$\mc{D}_{\alpha}(g,g)= \int_0^{+\infty} t^{2N-\alpha-1} \int_{\R^N} \abs{ h(t\cdot)* g }^2 dx dt \geq 0 $$
for a whatever nonnegative, radially symmetric $h \in C^{\infty}_c(\R^N)$ normalized in such a way that $\int_0^{+\infty} t^{N-\alpha-1} (h*h)(t) dt = C_{N,\alpha}$. 
\end{Proposition}

We state now the standard convergences for the nonlinear Choquard terms in the case of a subcritical growth. Hereafter we use $\lesssim$ to indicate less or equal up to a constant.
\begin{Proposition}
\label{prop_converg_generiche_nonloc}
 Assume \textnormal{(F1)} and \textnormal{(F2)}.
\begin{itemize}
\item Let $u_n \wto u$ in $H^s(\R^N)$. Then 
for any $\varphi \in H^s(\R^N)$ we have
$$\int_{\R^N} \big(I_{\alpha}*F(u_n)\big) f(u_n) \varphi \to \int_{\R^N} \big(I_{\alpha}*F(u)\big) f(u) \varphi. $$
\item Assume in addition \textnormal{(F3)}. 
Let $u_n \wto u$ in $H^s_r(\R^N)$. Then 
$$\int_{\R^N} \big(I_{\alpha}*F(u_n)\big) F(u_n) \to \int_{\R^N} \big(I_{\alpha}*F(u)\big) F(u) $$
and 
$$\int_{\R^N} \big(I_{\alpha}*F(u_n)\big) f(u_n) u_n \to \int_{\R^N} \big(I_{\alpha}*F(u)\big) f(u) u. $$
\end{itemize}
\end{Proposition}

\claim Proof.
Let $u_n \wto u$ in $H^s(\R^N)$, then (up to a subsequence) $u_n\to u$ in $L^p_{loc}(\R^N)$ for $p \in [1, 2^*_s)$, and a.$\,$e. pointwise.
By assumptions $F(u_n)$ is bounded in $L^{\frac{2N}{N+\alpha}}(\R^N)$ and $F(u_n) \to F(u)$ in $L^q_{loc}(\R^N)$ for $q \in [1, \frac{2N}{N+\alpha})$. 
Thus $F(u_n) \wto F(u)$ in $L^{\frac{2N}{N+\alpha}}(\R^N)$.
By some standard topological argument, the convergence holds for the whole sequence.

Moreover, by Proposition \ref{prop_HLS} we gain 
\begin{equation}\label{eq_dim_prima_forte}
I_{\alpha}*F(u_n) \wto I_{\alpha}*F(u) \quad \hbox{in $L^{\frac{2N}{N-\alpha}}(\R^N)$}.
\end{equation}
Let now $\varphi \in C^{\infty}_c(\R^N)$, and set $\Omega:= \supp(\varphi)$. Since $u_n \to u$ in $L^p(\Omega)$ for each $p \in [2, 2^*_s)$, we have (by the $L^p$-dominated convergence theorem) $f(u_n) \to f(u)$ in $L^p(\Omega)$ for each $p \in [1, \frac{2N}{\alpha+2s})$. Let $p \in (\frac{2N}{N+\alpha}, \frac{2N}{\alpha+2s})$ be whatever and let $q $ be such that $\frac{1}{p}+\frac{1}{q} = \frac{N+\alpha}{2N}$; since $\varphi \in L^q(\Omega)$ for such $q$, we have $f(u_n) \varphi \to f(u) \varphi$ in $L^{\frac{2N}{N+\alpha}}(\Omega)$, and actually in $L^{\frac{2N}{N+\alpha}}(\R^N)$. Thus, joining with \eqref{eq_dim_prima_forte},
$$ \int_{\R^N} \big(I_{\alpha}*F(u_n)\big) f(u_n) \varphi \to \int_{\R^N} \big(I_{\alpha}*F(u)\big) f(u) \varphi \quad \forall \varphi \in C^{\infty}_c(\R^N) .$$
We want to extend the relation to $\varphi \in H^s(\R^N)$. Indeed
\begin{align*}
\pabs{ \int_{\R^N} \big(I_{\alpha}*F(u_n)\big) f(u_n) \varphi} &\lesssim \norm{F(u_n)}_{\frac{2N}{N+\alpha}} \norm{f(u_n) \varphi}_{\frac{2N}{N+\alpha}} \\
&\lesssim \norm{|u_n|^{2^{\#}_{\alpha}-1}\varphi}_{\frac{2N}{N+\alpha}} + \norm{|u_n|^{2^*_{\alpha,s}-1} \varphi}_{\frac{2N}{N+\alpha}} 
\\
&\leq \norm{|u_n|^{\frac{\alpha}{N}}}_{\frac{2N}{\alpha}} \norm{\varphi}_{2} + \norm{|u_n|^{\frac{\alpha+2s}{N-2s}}}_{\frac{2N}{\alpha+2s}} \norm{\varphi}_{2^*_s} 
\\
&\leq \norm{u_n}_2^{2^{\#}_{\alpha}-1} \norm{\varphi}_{H^s} + \norm{u_n}_{2^*_s}^{2^*_{\alpha,s}-1} \norm{\varphi}_{H^s} 
 \lesssim \norm{\varphi}_{H^s},
\end{align*}
since $u_n$ are equibounded in $L^2(\R^N) \cap L^{2^*_s}(\R^N)$. Being the estimate uniform in $n$, we obtain the claim by density.

Assume now \hyperref[(F3)]{\textnormal{(F3)}}. Let $G(t):= (F(t))^{\frac{N+\alpha}{2N}}$. By the assumptions we have
$$\lim_{t \to 0} \frac{G(t)}{|t|^2} = \lim_{t \to 0} \left( \frac{F(t)}{|t|^{2^{\#}_{\alpha}}}\right)^{\frac{N+\alpha}{2N}}=0, \quad \quad \lim_{t \to \infty} \frac{G(t)}{|t|^{2^*_s}} = \lim_{t \to 0} \left( \frac{F(t)}{|t|^{2^{*}_{\alpha, s}}}\right)^{\frac{N+\alpha}{2N}}=0.$$
Thus, by \cite[Lemma 2.4]{ChWa0} we gain $G(u_n) \to G(u)$ in $L^1(\R^N)$, which means $F(u_n) \to F(u)$ in $L^{\frac{2N}{N+\alpha}}(\R^N)$. In particular, by Proposition \ref{prop_HLS} we obtain 
$$I_{\alpha}*F(u_n) \to I_{\alpha}*F(u) \quad \hbox{in $L^{\frac{2N}{N-\alpha}}(\R^N)$}.$$
Thus we get the first claim. Moreover, arguing as before we get $f(u_n)u_n \wto f(u) u$ in $L^{\frac{2N}{N+\alpha}}(\R^N)$, and this concludes the proof.
\QED

\medskip

\begin{Remark}\label{R:2.2}
We observe that, by substituting $F$ with $-F$, there is no loss of generality in assuming 
$$F(\sigma _0)>0 \quad \hbox{for some $\sigma_0 \neq 0$}$$
in \hyperref[(F4)]{\textnormal{(F4)}} ($\sigma_0$ can be chosen positive if, for example, \hyperref[(F5)]{\textnormal{(F5)}} holds) and 
$$\lim_{\sigma \to 0^+} \frac{F(\sigma )}{|\sigma|^{2^m_{\alpha,s}}} = + \infty$$
in \hyperref[(CF4)]{\textnormal{(CF4)}}. 
Thus, for the remaining part of the paper, we assume this positivity on the right-hand side of zero.
\end{Remark}


\section{Multidimensional annuli-shaped paths: even and odd nonlinearities}
\label{section:4..1}

We briefly denote by $q$ the lower-critical exponent $2^{\#}_{\alpha}$ and by $p$ the $L^2$-critical exponent $2^m_{\alpha,s}$, i.e.
$$q:=2^{\#}_{\alpha}=\frac{N+\alpha}{N}, \quad p:=2^m_{\alpha,s}= \frac{N+ \alpha+2s}{N}.$$

To avoid problems with the boundary of $\R_+$, we write from now on (see \cite{CGT2} for a different approach)
$$\mu \equiv e^{\lambda}\in (0, +\infty), \quad \lambda \in \R.$$
We also set
$$	\mc{D}(u) := \mc{D}_{\alpha}(F(u),F(u))=\int_{\R^N} (I_\alpha*F(u))F(u)\, dx.	$$
Using Proposition \ref{prop_HLS} and \hyperref[(F1)]{\textnormal{(F1)}}-\hyperref[(F2)]{\textnormal{(F2)}}, we notice that $\mc{D}$ is continuous on $L^2(\R^N) \cap L^{2^*_s}(\R^N)$, where $2^*_s= \frac{2N}{N-2s}$ is the Sobolev critical exponent, and thus continuous on $H^s_r(\R^N)$; notice that if we assume \hyperref[(CF2)]{\textnormal{(CF2)}}, then $\mc{D}$ is continuous also on $L^2(\R^N) \cap L^{2 + \frac{4s}{N+\alpha}}(\R^N)$.

To deal with the unconstrained problem, we further define the $C^1$-functional 
$\mc{J}:\R \times H^s_r(\R^N) \to \R$ by setting 
\begin{equation}\label{eq:2.3}
\mc{J} (\lambda, u):=\half \norm{(-\Delta)^{s/2} u}_2^2 -\half \mc{D}(u) + \frac{e^{\lambda}}{2} \|u\|_2^2, \quad (\lambda, u) \in \R\times H^s_r(\R^N).
\end{equation}
For a fixed $\lambda \in \R$, $u \in H^s_r(\R^N)$ is critical point of $\mc{J}(\lambda, \cdot)$ if and only if $u$ solves (weakly)
\begin{equation} \label{eq:2.4}
(- \Delta)^s u + e^{\lambda} u =(I_\alpha*F(u))f(u) \quad \hbox{in}\ \R^N.
\end{equation}

In this Section we study the geometry of
$$	u \in H^s_r(\R^N)\mapsto \mc{J}(\lambda,u) \in \R,	$$
for a fixed $\lambda\in\R$. We introduce a sequence of minimax values $a_n(\lambda)$, $n \in \N^*$: these values play important roles to find multiple solutions for the constrained problem (Theorem \ref{S:1.1_choq}) as well as for the unconstrained problem (Theorem \ref{S:1.3}).

For $n\in\N^*$ and $\lambda \in \R$ we introduce the set of paths
$$\Gamma_n(\lambda):=\big\{\gamma\in C(D_n, H^s_r(\R^N)) \mid \gamma \hbox{ odd}, \, \mc{J}(\lambda, \gamma_{|\partial D_n})<0 \big\},$$
where $D_n:=\{x \in \R^n \mid |x| \leq 1\}$, and the minimax values
$$a_n(\lambda):=\inf_{\gamma \in \Gamma_n(\lambda)} \sup_{\xi \in D_n} \mc{J}(\lambda, \gamma(\xi)).$$

For $n\geq 2$ the nonemptiness of $\Gamma_n(\lambda)$ has to be checked; for $n=1$ we refer to \cite[claim 1 of Proposition 2.1]{MS2}. 
Classically, in the local framework this fact was proved in \cite{BL2} by constructing inductively piecewise affine paths. This construction does not fit the nonlocality interaction given by the Choquard term, thus we need another approach. 

\begin{Proposition}\label{S:4.1}
Assume \hyperref[(F1)]{\textnormal{(F1)}}--\hyperref[(F4)]{\textnormal{(F4)}} and 
$F(\pm \sigma_0) \neq 0$.
Let $n \in \N^*$ and $\lambda \in \R$. 
Then $\Gamma_n(\lambda)\neq \emptyset$, thus $a_n(\lambda)$ is well defined. Moreover, $a_n(\lambda)>0$ and it is increasing with respect to $\lambda$ and $n$.
\end{Proposition}

\claim Proof.
Start observing that the polyhedron
 $$ \Sigma_n:=\left\{ t=(t_1,\dots,t_n) \mid \max_{i=1,\dots, n}\abs{t_i}=1\right\}
 $$
is homeomorphic to $\partial D_n$ (we pass from the $L^2$ to the $L^{\infty}$ norm). 
Let us fix $e_1, \dots, e_n \in C^{\infty}_c(\R^N)$, each of them between $0$ and $1$, radially symmetric, equal to one in some annulus $A_i$, and such that their supports are mutually disjoint. Then set $\gamma: \Sigma_n \to H^s_r(\R^N)$ by
\begin{equation}\label{eq_cammino_easy}
\gamma(t)(x):=\sigma_0 \sum_{i=1}^n t_i e_i(x)
\end{equation}
for every $t = (t_1, \dots, t_n) \in \Sigma_n$ and $x \in \R^N$. 
The map $\gamma$ is clearly odd and continuous. Moreover every $t \in \Sigma_n$ has at least a nontrivial component $|t_i|=1$, thus we have $F(\gamma(t)(x)) = F(\sigma _0 t_i e_i(x)) = F( \pm \sigma_0) \neq 0$ on $A_i$, hence $F(\gamma(t)) \nequiv 0$. By Proposition \ref{prop_positivity_Riesz} we have
$$\mc{D}(\gamma(t)) >0 \quad \hbox{ for each $t \in \Sigma_n$}.$$
Since
$\mc{D}\circ \gamma : \Sigma_n \to \R$ 
is continuous, and $\Sigma_n$ is compact, we obtain
$$ \min_{t \in \Sigma_n} \mc{D}(\gamma(t)) =: C>0,$$
i.e. 
$$
\mc{D}(\gamma(t)) \geq C >0 \quad \hbox{ for each $t \in \Sigma_n$}.
$$
Call moreover $M:= \max_{\Sigma_n} \norm{\gamma}_{H^s}^2 \in \R.$
By scaling, we obtain
\begin{align*}
\mc{J}_{\lambda}(\gamma(t)(\cdot/\theta))
&= \frac{\theta^{N-2s}}{2} \norm{(-\Delta)^{s/2} \gamma(t)}_2^2 + \frac{\theta^N e^{\lambda} }{2} \norm{\gamma(t)}_2^2 - \frac{\theta^{N+\alpha}}{2} \mc{D}(\gamma(t)) \\
&\leq \frac{\theta^{N-2s}}{2} M + \frac{\theta^N e^{\lambda} }{2} M - \frac{\theta^{N+\alpha}}{2} C < 0
\end{align*}
for some $\theta=\theta^* \gg 0$. Thus we consider $\tilde{\gamma}:=\gamma(\cdot)(\cdot/\theta^*) : \partial D_n \to H^s_r(\R^N)$. Finally we extend $\tilde{\gamma}$ to $D_n$ by
$$\tilde{\gamma}(\xi):= |\xi| \gamma\left(\frac{\xi}{|\xi|}\right)$$
for every $t \in D_n \setminus \{0\}$, and $\tilde{\gamma}(0):=0$. Therefore $\tilde{\gamma} \in \Gamma_n(\lambda) \neq \emptyset$.

What remains to prove is the monotonicity and positivity of $a_n(\lambda)$.
Since $D_n\subset D_{n+1}$, we may regard for $\gamma\in \Gamma_{n+1}(\lambda)$,
$$	\gamma_{|{D_n}}\in \Gamma_n(\lambda).	$$
Thus we have $a_n(\lambda)\leq a_{n+1}(\lambda)$. Since $\mc{J}(\lambda,u)$ is monotone
in $\lambda$, we also have the monotonicity with respect to $\lambda$.
The positivity of $a_1(\lambda)$ is instead essentially obtained in \cite{MS2}. 
Thus
$$ a_n(\lambda) \geq a_1(\lambda)>0.			 
\QED
$$

\bigskip

In the proof of Proposition \ref{S:4.1} we hardly relied on the positivity of the Riesz potential functional given in Proposition \ref{prop_positivity_Riesz}, to obtain the existence of path $\gamma: D_n \to H^s_r(\R^N)$ and a $C>0$ such that
\begin{equation}\label{eq_stima_path_C}
\mc{D}(\gamma(\xi)) \geq C >0 \quad \hbox{ for each $\xi \in \partial D_n$}.
\end{equation}
Anyway, no good information on $C$ are given by this approach.

A useful estimate in order to get infinitely many solutions for any $m>0$, when \hyperref[(CF4)]{\textnormal{(CF4)}} holds, is the one which relates $\mc{D}(\gamma(\theta \cdot))$ to $F(\theta \sigma_0)$ (see Lemma \ref{S:4.7} and Section \ref{sec_choquard_asympt_MP}), that is
\begin{equation}\label{eq_stima_path_F}
\mc{D}(\theta\gamma( \xi)) \geq C (F(\pm \theta \sigma_0))^2 \quad \hbox{ for each $\xi \in \partial D_n$ and $\theta \in [0,1]$}
\end{equation}
for some uniform $C>0$. When $F$ is nonnegative in a neighborhood of the origin (which is the case of $F$ even and \hyperref[(CF4)]{\textnormal{(CF4)}}), then one can easily build a suitable $\gamma$ which satisfies \eqref{eq_stima_path_F}.
\begin{Proposition}\label{prop_F_posit_intorn}
Assume \hyperref[(F1)]{\textnormal{(F1)}}--\hyperref[(F4)]{\textnormal{(F4)}}. 
Assume moreover that $F$ is nonnegative in some $[-\sigma_0, \sigma_0]$, $F(\pm \sigma_0) \neq 0$. Let $n \in \N^*$ and $\lambda \in \R$. Then the path $\gamma \in \Gamma_n(\lambda)$ defined in \eqref{eq_cammino_easy} satisfies \eqref{eq_stima_path_F}.
\end{Proposition}
\claim Proof.
Assume the notation of the proof of Proposition \ref{S:4.1}. For each $t\in \Sigma_n$, there exists $|t_{k}|=1$, thus, by exploiting that $\theta \sigma_0 t_i e_i(x) \in [-\sigma_0, \sigma_0]$ for each $x \in \R^N$, we obtain 
\begin{align*}
\mc{D}(\gamma(\theta t)) &= \int_{\R^N} \int_{\R^N} I_{\alpha}(x-y) F\left( \theta \sigma_0\sum_{i=1}^n t_i e_i (x)\right) F\left(\theta \sigma_0 \sum_{j=1}^n t_j e_j(x)\right) dx dy \\
& = \sum_{i=1}^n \sum_{j=1}^n \int_{\R^N} \int_{\R^N}I_{\alpha}(x-y) F\big(\theta \sigma_0 t_i e_i (x)\big) F\big( \theta \sigma_0 t_j e_j(x)\big) \, dx dy\\
& \geq \int_{\R^N} \int_{\R^N} I_{\alpha}(x-y) F\big(\theta \sigma_0 t_k e_k (x)\big) F\big( \theta \sigma_0 t_k e_k(y)\big) \, dx dy\\
& \geq \big( F( \pm \theta \sigma_0 )\big)^2 \int_{A_k} \int_{A_k} I_{\alpha}(x-y) \, dx dy
\end{align*}
which is the claim.
\QED

\bigskip

When $F$ is odd but not trivial around the origin (thus it cannot be here nonnegative) it seems not an easy task to build a $\gamma \in \Gamma_n(\lambda)$ satisfying \eqref{eq_stima_path_F}; indeed some estimate from below on 
$$\int_{\R^N} \int_{\R^N} I_{\alpha}(x-y) \frac{F\big(\theta \gamma(\xi)(x)\big) F\big(\theta \gamma(\xi)(y)\big)}{\big(F(\theta \sigma_0)\big)^2} dx dy$$
uniform for $\theta \to 0$ seems required; this is related to quotients of the type $\frac{F(\sigma h)}{F(\sigma )}$ with $\sigma \in (0,\sigma_0]$ and $h \in [0,1]$. 
This is essentially the meaning of condition \eqref{eq_cond_more_gen}.

\bigskip

To deal with this case we need a deep understanding of the Riesz potential on radial functions. 
We thus give now a different construction for a $\gamma \in \Gamma_n(\lambda)$: this procedure might be investigated also for more general Choquard-type equations, where different kernels (possibly sign-changing) appear.

We start by recalling a result contained in \cite[Theorem 1]{Thim0} (see also \cite[Lemma 6.3]{MMV} and references therein).

\begin{Theorem}[\cite{Thim0}]\label{S:4.2}
Let $\alpha \in (0, N)$ and $u\in L^1(\R^N) \cap L^{\infty}(\R^N) $ be radial. Then $I_{\alpha}*u$ is radial and
\begin{equation}
 \label{eq:4.1}
 (I_\alpha * u)(r)= \int_0^\infty F_\alpha\left(\frac{r}{\rho}\right) \rho^{\alpha-1}
 u(\rho)\, d \rho
\end{equation}
where $F_\alpha$ is positive and it satisfies, for some constants $C_{N,0}$, $C_{N,\infty}$, $C(N,\alpha)>0$, 
$$F_{\alpha}(\tau) \to C_{N,0}\ \hbox{ as $\tau \to 0$}, \quad 
F_{\alpha}(\tau) \tau^{N-\alpha}
 \to C_{N,\infty}\ \hbox{ as $\tau \to +\infty$}$$
and 
\begin{equation}
\frac{F_\alpha(\tau)}{G_\alpha(\tau)} \to C(N,\alpha) \quad \hbox{as}\ \tau\to 1, 
 \label{eq:4.2}
\end{equation}
with
\begin{equation}
 G_\alpha(s):=
\begin{cases}
 1 &\hbox{if $\alpha\in (1,N)$}, \cr
 \abs{\log\abs{\tau-1}} &\hbox{if $\alpha=1$}, \cr
 \abs{\tau-1}^{\alpha-1} &\hbox{if $\alpha\in (0,1)$}. \cr
\end{cases}
 \label{eq:4.3}
\end{equation}
\end{Theorem}

For a proof of Proposition \ref{S:4.1}, we prepare some notation and some estimates. We introduce the annuli
$$A(R,h):=\big\{x\in\R^N \mid \abs x\in [R-h,R+h]\big\}, \quad \chi(R,h;\cdot):= \chi_{A(R,h)}$$
for any $R\gg h >0$.

By applying Theorem \ref{S:4.2} to $u(|x|)\equiv\chi(1,h;\abs x)$ and arguing as in \cite[Lemma 1]{CGT4} (see also \cite{GalT}), we obtain 
\begin{align}
\lefteqn{\RRint I_\alpha(x-y)\chi(R,h;x)\chi(R,h;y)\, dxdy} \label{eq_quantity_to_bound} \\
 &= R^{N+\alpha}\RRint 
 I_\alpha(x-y)\chi\Big(1,\frac{h}{R};x\Big)
 \chi\Big(1,\frac{h}{R};y\Big)\, dxdy \notag \\
 &\sim 
\begin{cases}
 R^{N+\alpha}(\frac{h}{R})^2 &\hbox{if $\alpha\in (1,N)$}, \cr
 R^{N+1} (\frac{h}{R})^2\abs{\log\frac{h}{R}} &\hbox{if $\alpha=1$}, \cr
 R^{N+\alpha}(\frac{h}{R})^{1+\alpha} &\hbox{if $\alpha\in (0,1)$}.
 \end{cases} \notag 
\end{align}
For $R\geq 2$, we set the thickness of the annuli as
 $$ h_R :=
\begin{cases}
 R^{-\frac{N-2+\alpha}{2}} &\hbox{if $\alpha\in (1,N)$}, \cr 
 R^{-\frac{N-1}{2}}(\log R)^{-1/2} &\hbox{if $\alpha=1$}, \cr
 R^{-\frac{N-1}{1+\alpha}} &\hbox{if $\alpha\in (0,1)$,}
\end{cases}
 $$
so that a uniform bound for \eqref{eq_quantity_to_bound} is gained.

We show how to use these annuli to build a continuous odd map in $L^2(\R^N) \cap L^{2^*_s}(\R^N)$. 
By a regularization argument, we will obtain a refined path in $\Gamma_n(\lambda)$.

\medskip

\claim Proof of Proposition \ref{S:4.1} (refined).
In correspondence to $\sigma_0$, 
we will construct now a $\gamma\in \Gamma_n(\lambda)$; this path will moreover satisfy $\max_{\xi\in D_n} \norm{\gamma(\xi)}_{\infty} \leq \sigma_0.	$ 

\smallskip

\noindent 
\textbf{Step 1:} \emph{ Construction of an odd path in $L^r$.}
\\
For $n\geq 2$ we consider again the polyhedron
 $$ \Sigma_n=\big\{ t=(t_1,\dots,t_n) \mid \max_{i=1,\dots, n}\abs{t_i}=1\big\}.
 $$
For a large $R\gg 1$, which we will choose later, we define 
 $$ \gamma_R(t)(x):=\sum_{i=1}^n \hbox{sgn}(t_i) \chi\big(R^i, |t_i| h_{R^i};x\big):\,
\Sigma_n\to L^r(\R^N) $$
 where $r \in [1, +\infty]$. 
Here we regard $\chi(R^i,0;x)\equiv 0$.
We have 
\begin{align*}
 \calD(\sigma_0\gamma_R(t)) &= \sum_{i,j} 
 F(\hbox{sgn}(t_i)\sigma_0)F(\hbox{sgn}(t_j)\sigma_0) \cdot \cr
 &\quad\cdot \RRint I_\alpha(x-y)
 \chi(R^i,|t_i|h_{R^i}; x)
 \chi(R^i,|t_j|h_{R^i}; y)\, dxdy.
\end{align*}
We note that
\begin{itemize}
\item[(i)] For any $t=(t_1,\dots,t_n)\in\Sigma_n$, there exists at least one $t_k$ such that $\abs{t_k}=1$.
\item[(ii)] By \cite[Lemma 1]{CGT4}, 
 $$ F(\pm \sigma_0)^2 \RRint I_\alpha(x-y)
 \chi(R^k,h_{R^k};x)\chi(R^k,h_{R^k};y)
 \,dxdy \geq C_0.$$
\item[(iii)] By (i) and (ii),
 $$ \sum_{i=1}^n F(\pm \sigma_0)^2 \RRint I_\alpha(x-y)
 \chi(R^i,h_{R^i};x)\chi(R^i,h_{R^i};y)
 \,dxdy \geq C_0. $$
\item[(iv)] If $i\not=j$, by \cite[Lemma 2]{CGT4}, 
 $$ \RRint I_\alpha(x-y)\chi(R^i,h_{R^i};x)
 \chi(R^j,h_{R^j};x)\,dxdy \to 0
 \quad \hbox{as}\ R\to\infty. $$
\end{itemize}
By (i)--(iv), we have for sufficiently large $R\gg 1$,
\begin{equation}
 \calD(\sigma_0\gamma_R(t)) \geq C > 0 
\quad \hbox{for all}\ t\in\Sigma_n. \label{eq:4.6}
\end{equation}
In what follows we fix $R\gg 1$ so that \eqref{eq:4.6} holds.

\smallskip

\noindent
\textbf{Step 2:} \emph{Construction of an odd path in $H^s_r$.}
\\
For $0\leq h\ll R$ and $\epsilon>0$, we set
$$	\chi_\epsilon(R,h;x):=
\begin{cases}
1 &\text{if} \ x\in A(R,h), \\
1-\frac 1\epsilon \dista(x, A(R,h)), &\text{if}\ \dista(x,A(R,h))\in (0,\epsilon),\\
0 &\text{otherwise}.
 \end{cases}
$$
Here we regard
$	A(R,0)=\{ x\in\R^N|\, |x|=R\}.	$
We note that
\begin{align*}
&\chi_\epsilon(R,h;\cdot) \in H^s_r(\R^N) \ \text{for}\ \epsilon>0, \\
 	&\chi_\epsilon(R,h;\cdot)\to \chi(R,h;\cdot) \ \text{in}\ L^r(\R^N) \ \text{as}\ \epsilon\to 0 \ \text{for all}\ r\in [1,\infty), \\
&\supp\chi_\epsilon(R^i, h_{R^i};\cdot) \cap \supp\chi_\epsilon(R^j, h_{R^j};\cdot) =\emptyset \ \text{for}\ i\not=j \ \text{for $\epsilon$ small}.
\end{align*}
We set $\gamma_{\epsilon,R} \in C(\Sigma_n, H^s_r(\R^N))$ as
\begin{equation}\label{eq_def_gamma_epsR}
\gamma_{\epsilon,R}(t):=\sum_{i=1}^n \sgn(t_i)\chi_\epsilon(R^i, |t_i| h_{R^i}; \cdot), \quad t \in \Sigma_n. 
\end{equation}
 By \eqref{eq:4.6} and the continuity of $\mc{D}$ on $L^2(\R^N) \cap L^{2^*_s}(\R^N)$, we have for $\epsilon>0$ small
$$ \mc{D}(\sigma_0\gamma_{\epsilon,R}(t)) \geq C > 0 \quad \text{for all}\ t\in\Sigma_n. $$
Since
 $$ \mc{J}(\lambda, u(\cdot/\theta)) 
 = \half \theta^{N-2}\norm{(-\Delta)^{s/2} u}_2^2 +\frac{e^{\lambda}}{2}
 \theta^N\norm u_2^2 -\half\theta^{N+\alpha}\calD(u), $$
we have for large $\theta\gg 1$
 $$ \mc{J}(\lambda, \sigma_0\gamma_{\epsilon,R}(t)(\cdot/\theta))<0 \quad \hbox{for all}\ t\in\Sigma_n \approx \partial D_n. $$
Considering $D_n\equiv \{ ht \mid h\in [0,1], \, t\in\Sigma_n\}$ and extending $\sigma_0\gamma_{\eps,{R}}(t)(\cdot/\theta)$ to $D_n$ by
$$\widetilde \gamma(ht):=h{\sigma_0}\gamma_{\epsilon,R}(t)(\cdot/\theta),$$
finally we obtain a path $\widetilde\gamma\in \Gamma_n(\lambda)$. 
\QED

\medskip

\begin{Remark}
Even without assuming the positivity of $F$ (see Proposition \ref{prop_F_posit_intorn}), we notice that the construction of an odd map in $L^r$ gets easier when $F$ is an even function. 
Indeed there is no negative contribution given by the mixed interactions. We give only an outline of the proof, highlighting that in this case we do not need to use the fine Theorem \ref{S:4.2} given by \cite{Thim0}.

Define for every $i=1,\dots n$ and $h\in [-1,1]$ the annuli
$$A_i(h):= \big\{ x \in \R^N \mid |x| \in [2ni-|h|, 2ni+|h|]\big\},$$
and by $\chi_{A_i(h)}$ its indicator function. 
For every $t=(t_1, \dots, t_n) \in \Sigma_n$ we have that $A_1(t_1), \dots, A_n(t_n)$ are disjoint. Moreover, if $t_i=0$, then $\meas(A_i(t_i))=0$. 
Thus we define a continuous, odd map by
$$\gamma(t)(x):=\sum_{i=1}^n \sgn(t_i) \chi_{A_i(t_i)}(x): \, \Sigma_n\to L^2(\R^N)\cap L^{2^*_s}(\R^N). $$
Since $F$ is even, we obtain
\begin{align*}
 &\mc{D}(\sigma_0\gamma(t)) \\
 &= \sum_{i,j} \iint_{A_i(t_i)\times A_j(t_j)} 
I_{\alpha}(x-y) F(\sigma _0\sgn(t_i)\chi_{A_i(t_i)}(x)) F(\sigma _0\sgn(t_j)\chi_{A_j(t_j)}(y)) \,dxdy\\
 &= \big(F(\sigma _0)\big)^2 \sum_{i,j} \iint_{A_i(t_i)\times A_j(t_j)} I_{\alpha}(x-y) \,dxdy \geq\big(F(\sigma _0)\big)^2 C>0, 
\end{align*}
where $C$ does not depend on the specific $t$. The regularization to a $H^s_r$-path can be done as in the general case (or by mollification), as well as the extension to $D_n$.
\end{Remark}

We are ready now to show that $\gamma_{R,\epsilon}:\, \Sigma_n\to H^s_r(\R^N)$, defined in \eqref{eq_def_gamma_epsR}, has the desired property \eqref{eq_stima_path_F}.

\begin{Lemma} \label{S:4.7}
Assume \hyperref[(F1)]{\textnormal{(F1)}}--\hyperref[(F5)]{\textnormal{(F5)}}, and $F>0$ in some $(0,\delta_0]$.
If $F$ is odd, additionally assume \eqref{eq_cond_more_gen}. 
Then there exists a constant $A>0$ independent of $\sigma\in (0,\delta_0]$ and $t \in \Sigma_n$ such that
$$\mc{D}(\sigma\gamma_{R,\epsilon}(t)) \geq \frac{1}{2} \big( F(\sigma)\big)^2(A+o(1)) \quad \hbox{as}\ \epsilon\to 0. $$
Here $o(1)$ is a quantity which goes to $0$ as $\epsilon\to 0$ uniformly in $t\in\Sigma_n$ and $\sigma\in (0,\delta_0]$.
\end{Lemma}

\medskip

\claim Proof. 
We prove Lemma \ref{S:4.7} in two steps.

\smallskip

\noindent
\textbf{Step 1:} For $t\in\Sigma_n$, set
$$	a_{ij}(t):= \RRint I_\alpha(x-y)\chi(R^i,\abs{t_i}h_{R^i};x)
\chi(R^j,\abs{t_j}h_{R^j};y)\, dxdy.
$$
Then for sufficiently large $R>0$, we have
\begin{equation}\label{eq:4.9}	
A := \inf_{t\in\Sigma_n} \left(\sum_{i=1}^n a_{ii}(t)-\sum_{i\not=j}^n a_{ij}(t)\right)>0.
\end{equation}
This fact follows from \cite[Lemmas 1 and 2]{CGT4}. 
We fix $R\gg 1$ so that \eqref{eq:4.9} holds.

\smallskip

\noindent
\textbf{Step 2:} \emph{$\mc{D}(\sigma\gamma_{R,\epsilon}(t)) \geq\half F(\sigma )^2 A$ as $\epsilon\to 0$.} 
\\
We note that for $\epsilon>0$ small
$$	\supp\chi_\epsilon(R^i,\abs{t_i}h_{R^i}; \cdot) 
\cap \supp\chi_\epsilon(R^j,\abs{t_j}h_{R^j};\cdot)
=\emptyset \quad \hbox{for}\ i\not=j.	$$
Thus we have
\begin{align}
&\mc{D}(\sigma\gamma_{R,\epsilon}(t)) \nonumber \\
&= \sum_{i,j} \RRint I_\alpha(x-y)
F(\sigma \,\sgn(t_i)\chi_\epsilon(R^i,\abs{t_i}h_{R^i};x))
F(\sigma \,\sgn(t_j)\chi_\epsilon(R^j,\abs{t_j}h_{R^j};y))\, dxdy \nonumber \\
&=: \sum_{i,j} B_{ij}(\sigma,t).	\label{eq:4.10}
\end{align}
We consider cases $i=j$ and $i\not=j$ separately.

First we focus on the case $i=j$. For both even and odd $F$ we have
\begin{align}
&B_{ii}(\sigma,t) \nonumber \\
&= \RRint I_\alpha(x-y)
 		F(\sigma \,\sgn(t_i)\chi_\epsilon(R^i,\abs{t_i}h_{R^i};x))
F(\sigma \,\sgn(t_i)\chi_\epsilon(R^j,\abs{t_i}h_{R^i};y))\, dxdy \nonumber \\
&= \RRint I_\alpha(x-y)
 		F(\sigma \chi_\epsilon(R^i,\abs{t_i}h_{R^i};x))
F(\sigma \chi_\epsilon(R^j,\abs{t_i}h_{R^i};y))\, dxdy \nonumber \\
&\geq \RRint I_\alpha(x-y)
 		F(\sigma \chi(R^i,\abs{t_i}h_{R^i};x))
F(\sigma \chi(R^j,\abs{t_i}h_{R^i};y))\, dxdy \nonumber \\
&=F(\sigma )^2 a_{ii}(t),		\label{eq:4.11}
\end{align}
where we used the positivity of $F$ and the monotonicity of the integral. 
Next we consider the case $i\not=j$ for even $F$.
Since $F(\sigma )\geq 0$ for $\sigma\in [-\delta_0,\delta_0]$ we obtain
\begin{equation}\label{eq:4.12}
B_{ij}(\sigma,t) \geq 0 \quad \hbox{for all}\ t\in \Sigma_n.
\end{equation}
Finally we consider the case $i\not=j$ for odd $F$.
Since $\abs{F(\sigma )}=F(\abs \sigma)$ for $\sigma\in [-\delta_0,\delta_0]$
\begin{align}
&B_{ij}(\sigma,t) \nonumber \\
&= \RRint I_\alpha(x-y)
 		F(\sigma {\,\sgn(t_i)}\chi_\epsilon(R^i,\abs{t_i}h_{R^i};x))
F(\sigma {\,\sgn(t_j)}\chi_\epsilon(R^j,\abs{t_j}h_{R^j};y))\, dxdy \nonumber \\
&\geq - \RRint I_\alpha(x-y)
 		F(\sigma \chi_\epsilon(R^i,\abs{t_i}h_{R^i};x))
F(\sigma \chi_\epsilon(R^j,\abs{t_j}h_{R^j};y))\, dxdy. \label{eq:4.13}
\end{align}
Setting 
$$C_i(t,\epsilon):=\big\{x \mid \dista(x,A(R^i,\abs{t_i}h_{R^i}))\in (0,\epsilon)\big\}$$
we have
\begin{align*}
&\chi_\epsilon(R^i,\abs{t_i}h_{R^i};x) \in (0,1) \quad
\hbox{for}\ x\in C_i(t_i,\epsilon), \\
&\chi_\epsilon(R^i,\abs{t_i}h_{R^i};x) = \chi(R^i,\abs{t_i}h_{R^i};x) \quad
\hbox{for}\ x\not\in C_i(t_i,\epsilon),\\
&\meas(C_i(t_i,\epsilon)) \to 0 \quad \hbox{as}\ \epsilon\to 0, 
 {\hbox{ uniformly in $t \in \Sigma_n$}}.
\end{align*}
Thus for $r\in {[1,\infty)}$ {and $\sigma \in (0,\delta]$}
\begin{align}
\left\|\frac{1}{F(\sigma )}F(\sigma \chi_\epsilon(R^i,\abs{t_i}h_{R^i};\cdot))
-\chi(R^i,\abs{t_i}h_{R^i};\cdot)\right\|_r^r \nonumber 
&\leq \int_{C_i(t_i,\epsilon)} 
\pabs{\frac{1}{F(\sigma )}F(\sigma \chi_\epsilon(R^i,\abs{t_i}h_{R^i};x)) 
}^r\,dx 
\nonumber \\
&= \left(\max_{h\in [0,1]} \frac{\abs{F(\sigma h )}}{\abs{F(\sigma )}} 
\right)^r \meas(C_i(t_i,\epsilon))
\nonumber \\
&\to 0 \quad \hbox{as}\ \epsilon\to 0 \ \hbox{uniformly in}\ t\in\Sigma_n. \label{eq:4.14}
\end{align}
Here we use the fact that $\max_{h\in [0,1]} \frac{F(\sigma h)}{F(\sigma )} \leq C$, which follows from the local almost-monotonicity assumption in \hyperref[(CF4)]{\textnormal{(CF4)}}.
We note that \eqref{eq:4.14} implies, exploiting again \hyperref[(CF4)]{\textnormal{(CF4)}} 
\begin{align}
&\pabs{ \frac{1}{F(\sigma )^2}\RRint I_\alpha(x-y)
F(\sigma \chi_\epsilon(R^i,\abs{t_i}h_{R^i};x))
F(\sigma \chi_\epsilon(R^j,\abs{t_j}h_{R^j};y))\, dxdy -a_{ij}(t)} \nonumber \\
&\to 0 	\quad \hbox{as}\ \epsilon\to 0. \label{eq:4.15}
\end{align}
By \eqref{eq:4.13} and \eqref{eq:4.15},
\begin{equation}\label{eq:4.16}
B_{ij}(\sigma,t) \geq -F(\sigma )^2(a_{ij}(t)+o(1)) \quad \hbox{as}\ \epsilon\to 0.
\end{equation}
Thus, it follows from \eqref{eq:4.10}--\eqref{eq:4.12} and \eqref{eq:4.16} that
\begin{align*}
\mc{D}(\sigma\gamma_{R,\epsilon}(t)) 
&\geq F(\sigma )^2\left(\sum_{i=1}^na_{ii}(t)-\sum_{i\not=j}a_{ij}+o(1)\right) \\
&\geq \frac 12 F(\sigma )^2 (A + o(1))
\end{align*}
 This concludes the proof.
\QED


\section{Asymptotic analysis of mountain pass values}
\label{sec_choquard_asympt_MP}

We end this Section with some key estimates on the asymptotic behaviour of $a_n(\lambda)$ as $\lambda\to\pm\infty$.

\begin{Proposition}\label{S:4.6}
Assume \hyperref[(F1)]{\textnormal{(F1)}}--\hyperref[(F4)]{\textnormal{(F4)}} and let $n \in \N^*$.
\begin{itemize}
\item[(i)] If \hyperref[(CF3)]{\textnormal{(CF3)}} holds, then $\lim_{\lambda\to +\infty} \frac{a_n(\lambda)}{e^{\lambda}}=+\infty$.
\item[(ii)] If \hyperref[(CF4)]{\textnormal{(CF4)}} holds, then $\lim_{\lambda \to -\infty} \frac{a_n(\lambda)}{e^{\lambda}}=0$.
\end{itemize}
\end{Proposition}

Part (i) of Proposition \ref{S:4.6} can be found in \cite[Lemma 1]{CGT2}, once observed that $a_n(\mu)\geq a_1(\mu)$ for each $n \in \N^*$.

We deal now with the proof of (ii) of Proposition \ref{S:4.6}. 
We highlight that, when $F$ is even, the proof can be simplified (see \cite{CGT1} and Proposition \ref{prop_F_posit_intorn}).

The proof will be based on the key Lemma \ref{S:4.7}. 
We start noticing that, by \hyperref[(CF4)]{\textnormal{(CF4)}} and Remark \ref{R:2.2}, for some $\delta_0>0$ 
$$	F(\sigma )>0 \quad \hbox{for}\ \sigma\in (0,\delta_0],	$$
which implies
\begin{itemize}
\item[(i)] when $F$ is even, $F(\sigma )>0$ for all $\sigma\in [-\delta_0,\delta_0]\setminus\{ 0\}$;
\item[(ii)] when $F$ is odd, $F(\sigma )<0$ for all $\sigma\in [-\delta_0,0)$.
\end{itemize}
By \hyperref[(CF4)]{\textnormal{(CF4)}}, we also note that there exists $L_\sigma>0$ with $L_\sigma\to +\infty$ as
$\sigma\to 0^+$ such that
\begin{equation}\label{eq:4.8}	
F(\tau ) \geq L_\sigma \tau^p \quad \hbox{for all}\ \tau\in [0,\sigma].
\end{equation}

\claim Proof of (ii) of Proposition \ref{S:4.6}.
We write $\mu \equiv e^{\lambda}$ here. 
For $\sigma_0\in (0,\delta_0]$ and $\mu>0$, we consider the map
\begin{equation}\label{eq_dim_path_sigma}
\sigma t \in D_n \mapsto \sigma\sigma_0 \gamma_{R,\epsilon}(t)(\cdot/\mu^{-\frac 12}) \in H^s_r(\R^N) .
\end{equation}
We have by Lemma \ref{S:4.7} (since $\eps$ is fixed, we write $A$ instead of $A+o(1)$)
\begin{align*}
&\mu^{-1} \mc{J}(\mu,\sigma\sigma_0\gamma_{R,\epsilon}(t)(\cdot/\mu^{-\frac 12})) \\
&=\frac 12 \mu^{-\frac N2}(\sigma\sigma_0)^2 \norm{(-\Delta)^{s/2}\gamma_{R,\epsilon}(t)}_2^2
+\frac 12\mu^{-\frac N2 
}(\sigma\sigma_0)^2 \norm{\gamma_{R,\epsilon}(t)}_2^2
-\frac 12\mu^{-\frac N2 p}\mc{D}(\sigma\sigma_0\gamma_{R,\epsilon}(t)) \\
&\leq \frac 12\mu^{-\frac N2}(\sigma\sigma_0)^2 \norm{\gamma_{R,\epsilon}(t)}_{H^s}^2
-\frac 14 \mu^{-\frac N2 p} F(\sigma \sigma_0)^2A.
\end{align*}
Thus for $\mu$ small and $\sigma=1$
$$	\mc{J}(\mu,\sigma_0\gamma_{R,\epsilon}(t)(\cdot/\mu^{-\frac 12})) <0 
\quad \hbox{for}\ t\in\Sigma_n,	$$
which implies that \eqref{eq_dim_path_sigma} is a path in $\Gamma_n(\mu)$.
Moreover by \eqref{eq:4.8}
\begin{align*}
\mu^{-1}a_n(\mu) 
&\leq \max_{\sigma\in [0,1], \, t\in\Sigma_n} \mu^{-1}\mc{J}(\mu,\sigma\sigma_0\gamma_{R,\epsilon}(t)(\cdot/\mu^{-\frac 12}))
 \\
&\leq \max_{\sigma\in [0,1], \, t\in\Sigma_n} 
\frac 12\mu^{-\frac N2}(\sigma\sigma_0)^2 \norm{\gamma_{R,\epsilon}(t)}_{H^s}^2
-\frac 14 \mu^{-\frac N2 p} F(\sigma \sigma_0)^2A \\
&\leq \max_{\sigma\in [0,1], \, t\in\Sigma_n} 
\frac 12\mu^{-\frac N2}(\sigma\sigma_0)^2 \norm{\gamma_{R,\epsilon}(t)}_{H^s}^2
-\frac 14 L_{\sigma_0} (\mu^{-\frac N2} ( \sigma\sigma_0)^2)^p A \\
&\leq C_{\sigma_0},
\end{align*}
where
$$	C_{\sigma_0}:=\sup_{\tau\geq 0, \, t\in\Sigma_n} 
\left( \frac 12{\tau} \norm{\gamma_{R,\epsilon}(t)}_{H^s}^2
-\frac 14 L_{\sigma_0} A {\tau}^{p}\right) \in \R.	$$
Thus we have
$$	\limsup_{\mu\to 0^+} \mu^{-1} a_n(\mu) \leq C_{\sigma_0}.	$$
Since $C_{\sigma_0}\to 0$ as $\sigma_0\to 0$, we have (ii) of Proposition \ref{S:4.6}.
\QED


\section{Existence of multiple solutions} 
\label{sec_exist_sol}

\subsection{The Pohozaev mountain}
\label{section:4.3}

In this Section we start studying the Lagrangian formulation, applying the previous asymptotic estimates to a Pohozaev geometry.
We consider the functional 
$\mc{I}^m : \R \times H^s_r(\R^N) \to \R$ defined by
\begin{equation}\label{eq:2.2}
\mc{I}^m(\lambda, u):=\half \norm{(-\Delta)^{s/2} u}_2^2 -\half \mc{D}(u)+ \frac{e^{\lambda}}{2} \bigl( \|u\|_2^2 -m \bigr), \quad (\lambda, u) \in \R\times H^s_r(\R^N).
\end{equation}
It is immediate that, for any $(\lambda, u) \in \R\times H^s_r(\R^N)$,
$ \mc{I}^m(\lambda, u) = \mc{J}(\lambda, u) - \frac{e^{\lambda}}{2} m. $ 
Under some additional condition (see Section \ref{sec_Pohozaev} and \cite{DSS1}) every critical point satisfies the Pohozaev identity \eqref{eq:2.5}.
Inspired by this relation, we also introduce the Pohozaev functional $\mc{P}:\R \times H^s_r(\R^N) \to \R$ by setting 
\begin{equation}\label{eq:2.6}
\mc{P}(\lambda , u):= \frac{N-2s}{2}\|(-\Delta)^{s/2} u\|_2^2 -\frac{N+ \alpha}{2} \mc{D}(u) + \frac{N}{2} e^{\lambda} \|u\|_2^2 , \quad (\lambda, u) \in \R\times H^s_r(\R^N).
\end{equation}
We consider the action of $\G:=\Z_2$ on $\R^n$, $n \in \N^*$, and on $\R\times H^s_r(\R^N)$, given by
$$(\pm1,\xi) \in \G\times \R^n\mapsto \pm \xi \in \R^n,$$
$$ (\pm1,\lambda, u) \in \G\times \big(\R\times H^s_r(\R^N)\big) \mapsto (\lambda, \pm u)\in \R\times H^s_r(\R^N).$$
We notice that, under the assumption \hyperref[(F5)]{\textnormal{(F5)}}, $\mc{I}^m$, $\mc{J}$ and $\mc{P}$ are invariant under this action, i.e. they are even in $u$. 
In addition, we observe by the Principle of Symmetric Criticality of Palais \cite{Pal0} that every critical point of $\mc{I}^m$ restricted to $\R \times H^s_r(\R^N)$ is actually a critical point of $\mc{I}^m$ on the whole $\R \times H^s(\R^N)$.

\smallskip

Moreover we consider the \emph{Pohozaev set}
$$ \Omega :=\big\{(\lambda,u) \in \R \times H^s_r(\R^N) \mid \mc{P}(\lambda,u)>0\big\} \cup\big\{(\lambda,0) \mid \lambda \in \R \big\};$$
under the assumption \hyperref[(F5)]{\textnormal{(F5)}}, $\Omega$ is symmetric with respect to the axis $\{(\lambda,0) \mid \lambda \in \R \}$. 
We start showing the following property, due to the fact that $\mc{D}(u)=o(\norm u_{H^s}^2)$ as $u \to 0$.

\begin{Lemma}
We have
\begin{equation}\label{eq:4.17}	
\{ (\lambda,0) \mid \lambda\in \R\} \subset int(\Omega).
\end{equation}
\end{Lemma}

\claim Proof.
By
$$ |F(\sigma )| \lesssim |t|^q + |t|^p $$
where $q=\frac{N+\alpha}{N}$ and $p=\frac{N+ \alpha+2}{N}<2^*_s$.
 Thus 
$$ \|F(u)\|_{\frac{2N}{N+ \alpha}} \lesssim \| |u|^q\|_{\frac{2N}{N+ \alpha}} + \| |u|^p \|_{\frac{2N}{N+ \alpha}} 
= \| u\|_2^q + \| u \|_{\frac{2Np}{N+ \alpha}}^p . $$
Therefore by Proposition \ref{prop_HLS} and Young's inequality we have 
\begin{align*}
 \int_{\R^N} (I_\alpha*\abs{F(u)})\abs{F(u)}\,dx 
&\lesssim 	 \norm{F(u)}_{\frac{2N}{N+\alpha}}^2 
\lesssim \left( \norm{u}_2^q +\norm{u}_{\frac{2Np}{N+ \alpha}}^p\right)^2
 \\ 
&\lesssim \norm{u}_2^{2q} +\norm{u}_{\frac{2Np}{N+ \alpha}}^{2p} \leq \norm{u}_{H^s}^{2q} + \norm{u}_{H^s}^{2p}
\end{align*}
thus
$$ \mc{P}(\lambda,u) \lesssim \norm{u}_{H^s}^2 - \norm{u}_{H^s}^{2q} -\norm{u}_{H^s}^{2p} > 0$$
for $\norm{u}_{H^s}$ small, $u\neq 0$.
\QED

\bigskip

By \eqref{eq:4.17} we obtain
$$
\partial \Omega =\big\{(\lambda,u) \in \R \times H^s_r(\R^N) \mid \mc{P}(\lambda,u)=0, \ u \nequiv 0 \big\}
$$
and we call $\partial \Omega$ \emph{Pohozaev mountain}, on which $\mc{I}^m$ has height bounded from below (Proposition \ref{S:4.9}) and which separates zones where $\mc{I}^m$ is small (Proposition \ref{S:4.10}), as we see next. 
We observe that $\partial \Omega \neq \emptyset$, for instance by \cite[Theorems 1 and 3]{MS2}.

\begin{Proposition}\label{S:4.9}
Assume \hyperref[(F1)]{\textnormal{(F1)}}--\hyperref[(F4)]{\textnormal{(F4)}} and \hyperref[(F5)]{\textnormal{(F5)}}. We have the following properties.
\begin{itemize}
\item[(i)] $\mc{J}(\lambda,u)\geq 0$ for all $(\lambda,u)\in \Omega$.
\smallskip
\item[(ii)] $\mc{J}(\lambda,u)\geq a_1(\lambda)>0$ for all $(\lambda,u)\in \partial\Omega$. 
\item[(iii)] 
Assume \hyperref[(CF3)]{\textnormal{(CF3)}}. For any $m>0$, we set
$$E^m:= \inf_{(\lambda,u)\in\partial\Omega}\mc{I}^m(\lambda,u), \quad \text{ and } \quad B^m:=\inf_{\lambda \in \R} \left(a_1(\lambda) -\frac{e^\lambda}{2}m \right).$$
Then $E^m \geq B^m >-\infty$. 
\end{itemize}
\end{Proposition}

\claim Proof.
We notice that for all $(\lambda,u)\in \Omega$ 
$$ {\mathcal J}(\lambda,u) \geq {\mathcal J}(\lambda,u) - \frac{{\mathcal P}(\lambda,u)}{N + \alpha} 
= \frac{\alpha +2}{2(N + \alpha)} \norm{(-\Delta)^{s/2} u}_2^2+ \frac{\alpha}{2(N + \alpha)} e^{\lambda}\norm{u}_2^2 \geq 0 $$
and thus $(i)$ follows. 
As regards $(iii)$: the fact that $E^m \geq B^m$ is a direct consequence of $(ii)$, while the fact that $B^m >-\infty$ comes from Proposition \ref{S:4.6} (i).

Focus now on point $(ii)$: this comes from the fact that for each $\lambda$ the mountain pass level $a_1(\lambda)$ coincides with the ground state energy level (see \cite[Section 4.2]{MS1} and \cite[Proposition 3.3]{CGT3} for details).
Differently, without exploiting the existence result for the unconstrained problem, we argue as follows (see \cite[Remark 3.4.3]{GalT}): 
let $\gamma \in \Gamma(\lambda)$; by definition of $\Gamma_1(\lambda)$ and by point (i) there exists $t^*$ such that $\gamma(t^*) \in \partial \Omega$ and $\gamma(t^*) \neq 0$, thus $\mc{P}(\lambda, \gamma(t^*))=0$. This means that
$$\mc{J}(\lambda, \gamma(t^*)) =\frac{\alpha +2s}{2(N + \alpha)} \norm{(-\Delta)^{s/2} \gamma(t^*)}_2^2+ \frac{\alpha \mu}{2(N + \alpha)} \norm{\gamma(t^*)}_2^2 \simeq \norm{\gamma(t^*)}_{H^s}^2 $$
thus
$$a(\lambda) \gtrsim \inf_{u \in (\partial \Omega)_{\lambda}} \norm{u}_{H^s}^2$$
where $(\partial \Omega)_{\lambda} :=
\big\{ u \in H_r^s (\R^N) \setminus \{ 0 \} \mid \mc{P}(\lambda,u)=0 \big\}$. 
Since, by \eqref{eq:4.17}, $ (\partial \Omega)_{\lambda}$ is far from the line $(\lambda,0)$, we obtain that the right-hand side is strictly positive, which is the claim.
\QED

\bigskip

 From now on we assume \hyperref[(CF3)]{\textnormal{(CF3)}} to give sense to the quantity $B^m$. In view of Proposition \ref{S:4.9} (iii), 
we set for $m>0$ and $n\in\N^*$
\begin{align*}
\Gamma_n^m:=\big\{\Theta \in C(D_n, \R\times H^s_r(\R^N)) \mid \; & \text{$\Theta$ is $\G$-equivariant,} \ \mc{I}^m(\Theta(0)) \leq B^m-1 ,\\
& \Theta|_{\partial D_n}\notin \Omega, \ \mc{I}^m(\Theta|_{\partial D_n})\leq B^m-1 \big\}
\end{align*}
and
$$b_n^m := \inf_{\Theta \in \Gamma_n^m} \sup_{\xi \in D_n} \mc{I}(\Theta(\xi));$$
we point out that asking $\Theta=(\Theta_1, \Theta_2) \in \Gamma^m_n$ to be $\G$-equivariant means that $\Theta_1$ is even and $\Theta_2$ is odd, and in particular $\Theta_2(0)=0$ which implies $\Theta(0) \in \Omega$.
Arguing as in \cite{CGT1,CGT4} we obtain the following.

\begin{Proposition}\label{S:4.10}
Assume \hyperref[(F1)]{\textnormal{(F1)}}-\hyperref[(F2)]{\textnormal{(F2)}}-\hyperref[(CF3)]{\textnormal{(CF3)}}-\hyperref[(F4)]{\textnormal{(F4)}}-\hyperref[(F5)]{\textnormal{(F5)}}. We have the following properties.
\begin{itemize}
\item[(i)] For any $m>0$ and $n \in \N^*$, we have $\Gamma_n^m\neq \emptyset$ and
\begin{equation*}\label{eq:4.18}
b_n^m\leq a_n(\lambda) - e^{\lambda} \frac{m}{2},
\end{equation*}
for each $\lambda \in \R$. Moreover, $b_n^m$ increases with respect to $n$.
\item[(ii)] For any $k\in\N^*$ there exists $m_k\geq 0$, namely given by
\begin{equation*}\label{eq:4.19}
m_k:= 2\inf_{\lambda\in\R}\frac{a_k(\lambda)}{e^{\lambda}},
\end{equation*}
such that for $m>m_k$
$$	b_n^m < 0 \quad \text{for}\ n=1,2,\dots, k.	$$
Moreover, $m_k$ is increasing with respect to $k$.
\item[(iii)] If \hyperref[(CF4)]{\textnormal{(CF4)}} holds, then $m_k=0$ for each $k \in \N^*$. That is, for each $m>0$ we have
$$	b_n^m < 0 \quad \text{for all}\ n\in\N^*.	$$
\end{itemize}
\end{Proposition}

\begin{Corollary}
For any $m>0$, we have 
$$B^m= E^m = b^m_1,$$
i.e. the first minimax value $b_1^m$ equals the Pohozaev minimum $E^m$ on the product space.
\end{Corollary}


\subsection{The Palais-Smale-Pohozaev condition} 
\label{sec_PSP_choquard}

For every $b \in \R$ we set
$$ K_b^{P,m} := \big\{ (\lambda,u)\in \R\times H^s_r(\R^N) \mid \mc{I}^m(\lambda, u)=b,\,
 D \mc{I}^m(\lambda, u)=0, \,
 \mc{P}(\lambda, u)=0 \big\}.$$
 We notice that, assuming \hyperref[(F5)]{\textnormal{(F5)}}, $K_b^{P,m}$ is invariant under the $\G$-action.

As already observed, under \hyperref[(F1)]{\textnormal{(F1)}}-\hyperref[(F2)]{\textnormal{(F2)}} it is not known if $\partial_u\mc{I}^m(\lambda, u)=0$ implies $\mc{P}(\lambda, u)=0$, thus it seems hard to recognize the compactness of $K^{m}_b := \big\{ (\lambda,u)\in \R\times H^s_r(\R^N) \mid \mc{I}^m(\lambda, u)=b,\, D \mc{I}^m(\lambda, u)=0\big\} $ (through the standard Palais-Smale condition) without the additional condition $\mc{P}(\lambda, u)=0$.

Inspired by \cite{HIT,HT0,IT0}, we make use of the Palais-Smale-Pohozaev condition, a weaker compactness condition that takes into account the scaling properties of $\mc{I}^m$ through the Pohozaev functional $\mc{P}$. 
Through this tool we will show that $K_b^{P,m}$ is compact when $b<0$.

We thus have the following result, which can be found in \cite[Theorem 3]{CGT2}.

\begin{Proposition}\label{S:3.2}
Assume \hyperref[(F1)]{\textnormal{(F1)}}-\hyperref[(CF2)]{\textnormal{(CF2)}}-\hyperref[(CF3)]{\textnormal{(CF3)}} and let $b <0$. Then $\mc{I}^m$ satisfies the the \emph{Palais-Smale-Pohozaev condition} at level $b$ (shortly the $(PSP)_b$ condition), that is
 every $(\lambda_n, u_n)_n \subset \R \times H^s_r(\R^N)$ satisfying 
\begin{equation}\label{eq:3.1}
\mc{I}^m (\lambda_n, u_n) \to b,
\end{equation}
\begin{equation}\label{eq:3.2}	
\partial_{\lambda} 	\mc{I}^m(\lambda_n, u_n) \to 0, 
\end{equation}
\begin{equation}\label{eq:3.3}
\norm{\partial_u 	\mc{I}^m(\lambda_n, u_n)}_{(H^s_r(\R^N))^*} \to 0,
\end{equation}
\begin{equation}\label{eq:3.4}
\mc{P}(\lambda_n, u_n) \to 0.
\end{equation}
has a strongly convergent subsequence in $\R \times H^s_r(\R^N)$.

As a consequence, $K_b^{P,m} \cap (\R \times \{0\}) = \emptyset$ and $K_b^{P,m}$ is compact.
\end{Proposition}

We emphasize that $K^{P,m}_0$ is not compact, since we can consider an unbounded sequence $(\lambda_n, 0)$ with $\lambda_n \to - \infty$.

\medskip

\claim Proof of Theorem \ref{S:1.1_choq}.
Once obtained the $(PSP)$ condition, one can argue as in \cite{CGT1,CGT4,HT0}. 
We construct a deformation flow for $\mc{I}^m$ through a deformation of an augmented functional
 $$ \mc{H}^m(\theta,\lambda,u) :=\mc{I}^m(\lambda,u(e^{-\theta}\cdot))
 $$
on a Hilbert manifold $M:=\R\times\R\times H_r^s(\R~N)$ with a suitable metric.
Applying genus arguments we can find multiple solutions.
We refer to \cite{CGT1,CGT4,HT0} for details.
\QED

\medskip

\claim Proof of Theorem \ref{S:1.3}.
Theorem \ref{S:1.3} can obtained as a byproduct of the techniques developed for Theorem \ref{S:1.1_choq}. See \cite[Section 5]{CGT4} for details.
\QED


\section{The Pohozaev identity} 
\label{sec_Pohozaev}

In \cite[equation (6.1)]{DSS1}, in presence of power nonlinearities, it is proved that every weak solution $u$ is $C^2$ and thus satisfies the Pohozaev identity \eqref{eq:2.5}.

Here we want to extend the identity to more general nonlinearities and to more general solutions $u\in C^1$.
We consider the following assumption:
\begin{itemize}
\item[\textnormal{(F6)}] \label{(F6)}
$f\in C^{0,\sigma}_{loc}(\R)$ for some $\sigma \in (0,1]$. 
\end{itemize}
From \cite[Theorem 1.2 and Proposition 4.9]{CGT3} and \cite[Section 3.1]{CG1} 
(see also \cite[Section 4.4.5]{GalT}) we have the following regularity results.

\begin{Proposition}[\cite{CGT3,CG1,GalT}]
\label{prop_u_L1}
Assume that \hyperref[(F1)]{\textnormal{(F1)}}-\hyperref[(F2)]{\textnormal{(F2)}} hold. Let $u\in H^s(\R^N)$ be a weak solution of \eqref{eq_Choquard_genericaF}. Then
\begin{itemize}
\item[i)] if $u$ is nonnegative or $\limsup_{\sigma\to 0} \frac{|f(\sigma)|}{|\sigma|}<\infty$, then $u \in L^{\infty}(\R^N)$.
\end{itemize}
 Let now $u\in H^s(\R^N)\cap L^{\infty}(\R^N)$ be a weak solution of \eqref{eq_Choquard_genericaF}. Then
\begin{itemize}
\item[ii)] $u \in L^1(\R^N)\cap H^{2s}(\R^N)$; in particular, $F(u) \in L^1(\R^N) \cap L^{\infty}(\R^N)$;
\item[iii)] if $s\in (0,\frac{1}{2}]$, then $u \in C^{0,\gamma}(\R^N)$ for any $\gamma \in (0, 2s]\cap(0,1)$,
\item[iv)] if $s\in (\tfrac{1}{2}, 1)$, then $u \in C^{1,\gamma}(\R^N)$ for any $\gamma \in (0, 2s-1)$.
\end{itemize}
Assume \hyperref[(F6)]{\textnormal{(F6)}} in addition. Then
\begin{itemize}
\item[v)] if $s \in (0, \frac{1}{2})$ and $\omega:=\min\{2s \sigma, \alpha\} \leq 1-2s$, then $u \in C^{0,\gamma}(\R^N)$ for every $\gamma \in (0,\omega+2s]\cap (0,1)$;
\item[vi)] if $s \in [\frac{1}{4}, \frac{1}{2})$ and $\omega:=\min\{2s \sigma, \alpha\} > 1-2s$, then $u \in C^{1,\gamma}(\R^N)$ for every $\gamma \in (0,\omega+2s-1]\cap (0,1)$;
\item[vii)] if $s \in [\frac{1}{2}, 1)$ and $\omega:=\min\{\sigma, \alpha\}\leq 2-2s$, then $u \in C^{1,\gamma}(\R^N)$ for any $\gamma \in (0, \omega+2s-1] \cap (0,1)$;
\item[viii)] if $s \in (\frac{1}{2}, 1)$ and $\omega:=\min\{\sigma, \alpha\} > 2-2s$ then $u \in C^{2,\omega+2s-2}(\R^N)$;
\end{itemize}
and moreover
\begin{itemize}
\item[ix)] if $\alpha<2$ and $\sigma > 1-2s$,
then $u\in C^{1, \gamma}(\R^N)$ for every $\gamma \in (0,1)$.
\end{itemize}
Finally assume $f \in C^1(\R)$. Then
\begin{itemize}
\item[x)] if $s \in (\frac{1}{2}, 1)$, then $u \in C^{2,\gamma}(\R^N)$ for every $\gamma\in (0,2s-1)$.
\end{itemize}
\end{Proposition}

In particular, we have the following.
\begin{Corollary}
\label{corol_reg_pohz}
Assume that \hyperref[(F1)]{\textnormal{(F1)}}-\hyperref[(F2)]{\textnormal{(F2)}} hold. 
Let $u\in H^s(\R^N)\cap L^{\infty}(\R^N)$
be a weak solution of \eqref{eq_Choquard_genericaF}. Assume in addition one of the following
\begin{itemize}
\item $s \in (\frac{1}{2}, 1)$,
\item $s=\frac{1}{2}$ and \hyperref[(F6)]{\textnormal{(F6)}},
\item $s \in [\frac{1}{4}, \frac{1}{2})$, $\alpha \in ( 1-2s, N)$ and \hyperref[(F6)]{\textnormal{(F6)}} with $\sigma \in (\frac{1-2s}{2s}, 1]$,
\item $s \in (0, \frac{1}{2})$, $\alpha \in (0,2)$ and \hyperref[(F6)]{\textnormal{(F6)}} with $\sigma \in (1-2s, 1]$.
\end{itemize}
Then $u \in C^{1,\gamma}(\R^N)$ for some $\gamma \in (0,1)$. If $s \in (\frac{1}{2}, 1)$ and \hyperref[(F6)]{\textnormal{(F6)}} holds too, then we can choose $\gamma \in (2s-1,1)$.
\end{Corollary}

\smallskip

Under a pointwise well posedness of the fractional Laplacian, and the existence of a weak gradient, we obtain the following integration by parts rule, inspired by \cite[Lemma 4.2]{DFW1}.

\begin{Proposition}
\label{prop_int_kern_div}
Let $s\in (0,1)$. 
Let $u \in D^{s,2}(\R^N) \cap C^{\gamma}_{loc}(\R^N) \cap Lip_{loc}(\R^N)$ for some $\gamma >2s$, and assume \eqref{eq_spazio_gener_s}.
Let moreover $X \in C^1_c (\R^N,\R^N)$ be a vector field, and define, for $x, y\in \R^N$, $x \neq y$, 
$$\mc{K}_{X}^s(x,y):= \frac{\big(\dive (X)\big)(x) + \big(\dive (X)\big)(y)}{2}- \frac{N+2s}{2}\frac{(X(x)-X(y))\cdot(x-y)}{|x-y|^2}$$
the \emph{fractional divergence kernel} related to $X$. 
Then it holds
$$\frac{C_{N,s}}{2}\int_{\R^N} \int_{\R^N} \frac{|u(x)-u(y)|^2}{|x-y|^{N+2s}} \mc{K}_{X}^s(x,y) dx dy = -
 \int_{\R^N} (-\Delta)^s u \, (\nabla u \cdot X ) \, dx;$$
noticed that the left-hand side is the weighted Gagliardo seminorm with weight $\mc{K}_{X}^s$, set 
$$\mc{G}_u^s(x,y):= \frac{C_{N,s}}{2} \frac{|u(x)-u(y)|^2}{|x-y|^{N+2s}}$$
we can write
$$ (\mc{G}_u^s, \mc{K}_X^s)_{L^2(\R^{2N})} =
- \big((-\Delta)^{s} u \, \nabla u,X\big)_{L^2(\R^N)}.$$
\end{Proposition}

\claim Proof. 
For the proof, we follow the lines of \cite{Dje0} (see also \cite{DeNDj}).
We start noticing that, being $u \in D^{s,2}(\R^N)$, by the assumptions we have
$$\mc{G}_u^s \in L^1(\R^{2N}), \quad \mc{K}^s_X \in L^{\infty}(\R^{2N})$$
so that the product is summable. By Dominated Convergence theorem, the symmetry of the kernel, and the Fubini theorem, we obtain
\begin{align*}
\MoveEqLeft 
\frac{2}{C_{N,s}}(\mc{G}_u^s, \mc{K}_X^s)_{L^2(\R^{2N})} \\
&= \lim_{\eps \to 0} \iint_{|x-y|>\eps} \frac{|u(x)-u(y)|^2}{|x-y|^{N+2s}} \mc{K}_{X}^s(x,y) dx dy \\
&= \lim_{\eps \to 0} \iint_{|x-y|>\eps}\frac{|u(x)-u(y)|^2}{|x-y|^{N+2s}} \left( \dive(X)(x) - (N+2s) \frac{(x-y)\cdot X(x)}{|x-y|^2} \right)dxdy \\
&= \lim_{\eps \to 0} \int_{\R^N} \left( \int_{\R^N \setminus B_{\eps}(y)}\frac{|u(x)-u(y)|^2}{|x-y|^{N+2s}} \left( \dive(X)(x)- (N+2s)\frac{(x-y)\cdot X(x)}{|x-y|^2} \right) dx \right) dy.
\end{align*}
Exploiting that, for $x\neq y$, $\nabla_x \frac{1}{|x-y|^{N+2s}} =- (N+2s) \frac{x-y}{|x-y|^{N+2s+2}}$, and the Divergence theorem (possible because $\frac{X}{|\cdot-y|^{N+2s}} \in C^1_c (\overline{\R^N \setminus B_{\eps}(y)})$ and $u \in Lip_{loc}(\R^N) \subset W^{1,\infty}(\supp(X))$, see \cite[Theorem 4.6]{EG0})
\begin{align*}
\MoveEqLeft 
\frac{2}{C_{N,s}} (\mc{G}_u^s, \mc{K}_X^s)_{L^2(\R^{2N})} \\ 
=& \lim_{\eps \to 0} \int_{\R^N} \left(\int_{\R^N \setminus B_{\eps}(y)}|u(x)-u(y)|^2 \left( \frac{\dive(X)(x)}{|x-y|^{N+2s}} - (N+2s)\frac{(x-y)\cdot X(x)}{|x-y|^{N+2s+2}} \right) dx \right) dy \\
=&\lim_{\eps \to 0} \int_{\R^N} \left(\int_{\R^N \setminus B_{\eps}(y)}|u(x)-u(y)|^2 \dive_x \left( \frac{X}{|x-y|^{N+2s}}\right)(x) dx \right) dy \\
=&- 2 \lim_{\eps \to 0} \int_{\R^N} \left(\int_{\R^N \setminus B_{\eps}(y)}(u(x)-u(y)) \nabla u(x) \cdot \frac{X(x)}{|x-y|^{N+2s}} dx \right) dy + \\
& +\lim_{\eps \to 0} \int_{\R^N} \left(\int_{\partial B_{\eps}(y)} |u(x)-u(y)|^2 \frac{X(x)}{|x-y|^{N+2s}} \cdot \frac{x-y}{|x-y|} d\sigma(x) \right) dy \\
=&- 2\lim_{\eps \to 0} \int_{\R^N} \left(\int_{\R^N \setminus B_{\eps}(y)} \frac{u(x)-u(y)}{|x-y|^{N+2s}} \nabla u(x) \cdot X(x) dx \right) dy + \\
& + \lim_{\eps \to 0} \frac{1}{\eps^{N+2s+1}} \int_{\R^N} \left(\int_{\partial B_{\eps}(y)} |u(x)-u(y)|^2 X(x) \cdot (x-y) d\sigma(x) \right) dy.
\\
=:& \, -2 \lim_{\eps \to 0} I_{\eps} + \lim_{\eps \to 0} E_{\eps};
\end{align*}
here we split the limits since we will prove the existence of both. 

For the first integral, we notice that $x \mapsto \int_{\R^N \setminus B_{\eps}(x)} \frac{|u(x)-u(y)|}{|x-y|^{N+2s}} |\nabla u(x) \cdot X(x)| dy \leq C_{\eps} \norm{u}_{\infty} |\nabla u (x)\cdot X(x)| \in L^1(\R^N)$ so that we can apply Fubini theorem, 
then we perform a symmetrization substitution and apply again Fubini theorem, and finally Dominated Convergence theorem (since $y \mapsto \frac{ 2u(x)-u(x+y)-u(x-y)}{|y|^{N+2s}} \in L^1(\R^N)$ by Proposition \ref{prop_well_posed}), obtaining 
\begin{align*}
C_{N,s} \lim_{\eps \to 0} I_{\eps} &= C_{N,s} \lim_{\eps \to 0} \iint_{|x-y|>\eps} \frac{u(x)-u(y)}{|x-y|^{N+2s}} \nabla u(x) \cdot X(x) dx dy \\
&= \frac{C_{N,s}}{2}\lim_{\eps \to 0} \iint_{|y|>\eps} \frac{2u(x)-u(x+y)-u(x-y)}{|y|^{N+2s}} \nabla u(x) \cdot X(x) dx dy \\
&= \lim_{\eps \to 0} \int_{\R^N} \nabla u(x) \cdot X(x) \left(\frac{C_{N,s}}{2}\int_{\R^N \setminus B_{\eps}(0)} \frac{ 2u(x)-u(x+y)-u(x-y)}{|y|^{N+2s}} dy \right) dx \\
&=\int_{\R^N} \nabla u \cdot X (-\Delta)^s u.
\end{align*}

For the second integral, notice that the set $\{(x,y) \in \R^{2N} \mid x \in \supp(X), \, |x-y|=\eps\}$ is bounded.
Thus the integrand (being bounded) is summable, which allows us to implement the Fubini theorem and obtain, by exploiting also a symmetrization argument, 
\begin{align*}
(N+2s)E_{\eps} &= \frac{1}{\eps^{N+2s+1}} \iint_{|x-y|=\eps} |u(x)-u(y)|^2 X(x) \cdot (x-y) d\sigma(x) \times dy \\
&= \frac{1}{2\eps^{N+2s+1}} \iint_{|x-y|=\eps}|u(x)-u(y)|^2 (X(x)-X(y)) \cdot (x-y) d\sigma(x) \times dy .
\end{align*}
If $\supp(X)\subset B_R(0)$, then out of the set
$$A_{R,\eps}:=\{(x,y) \in B_R(0)\times B_R(0) \mid |x-y|=\eps\}$$
the integrand is null. 
Thus, being $u \in Lip_{loc}(\R^N)$ (actually it is sufficient $u \in C^{0,\theta}_{loc}(\R^N)$ for some $\theta>s$) and $X \in Lip(\R^N, \R^N)$, we get
\begin{align*}
E_{\eps} &\lesssim \frac{1}{\eps^{N+2s+1}} \iint_{A_{R,\eps}} |x-y|^4 d\sigma(x) \times dy \\
&=\eps^{-N-2s +3} m_{2N-1}(A_{R,\eps}).
\end{align*}
Observed that $m_{2N-1}(A_{R,\eps}) \lesssim m_N(B_R) m_{N-1}(\partial B_{\eps}) \sim \eps^{N-1}$, we obtain $E_{\eps} \lesssim \eps^{-2s+2} \to 0$. Joining the pieces, we reach the claim.
\QED

\medskip

\begin{Corollary}
In the assumptions of Proposition \ref{prop_int_kern_div}, let $G \in C^1(\R^N)$ with $G(u) \in L^1(\R^N)$ and 
$$(-\Delta)^s u = g(u) \quad \hbox{ in $\R^N$}$$
in the pointwise sense, where $G'=g$. Then 
\begin{align*} 
\frac{C_{N,s}}{2} \int_{\R^N} \int_{\R^N} \frac{|u(x)-u(y)|^2}{|x-y|^{N+2s}} \mc{K}_{X}^s(x,y) dx dy &= -
 \int_{\R^N} \nabla G(u) \cdot X \, dx \\
 &= \int_{\R^N} G(u) \, \dive(X) dx, 
 \end{align*}
 i.e.
 $$ (\mc{G}_u^s, \mc{K}_X^s)_{L^2(\R^{2N})} =
(G(u), \dive(X))_{L^2(\R^N)}.$$
\end{Corollary}

\smallskip

We deal now with the Riesz kernel, right-hand side of the equation.

\begin{Proposition}
\label{prop_int_riesz_div}
Let $\alpha \in (0,N)$ and $H \in Lip_{loc}(\R^N)\cap L^{\infty}(\R^N)$ be such that 
$$ (I_{\alpha}*|H|)|H| \in L^1(\R^N), \quad (I_{\alpha}*|H|)|\nabla H| \in L^1_{loc}(\R^N).$$
Let moreover $X \in C^1_c (\R^N,\R^N)$ be a vector field and set, for $x, y\in \R^N$, $x \neq y$, 
$$\mc{K}_{X}^{-\frac{\alpha}{2}}(x,y):= \frac{\big(\dive (X)\big)(x) + \big(\dive (X)\big)(y)}{2}- \frac{N-\alpha}{2} \frac{(X(x)-X(y))\cdot(x-y)}{|x-y|^2}.$$
Then 
\begin{align*} 
\int_{\R^N} \int_{\R^N} I_{\alpha}(x-y) H(x) H(y) \mc{K}_{X}^{-\frac{\alpha}{2}}(x,y) dx dy &= -
 \int_{\R^N}\big (I_{\alpha}*H\big) \nabla H \cdot X \, dx, 
 \end{align*}
 i.e. set
$$\mc{R}_H^{\alpha}(x,y):= I_{\alpha}(x-y) H(x) H(y)$$
we have
$$ (\mc{R}_{H}^{\alpha}, \mc{K}_{X}^{-\frac{\alpha}{2}})_{L^2(\R^{2N})} = - \big((I_{\alpha}*H) \nabla H,X\big)_{L^2(\R^N)}.
$$
\end{Proposition}

\claim Proof. 
We proceed as in the proof of Proposition \ref{prop_int_kern_div}.
We start noticing that 
$$\mc{R}_{H}^{\alpha} \in L^1(\R^{2N}), \quad \mc{K}_{X}^{-\frac{\alpha}{2}} \in L^{\infty}(\R^{2N})$$
by the assumptions, so that the product is summable. Thus
\begin{align*}
\MoveEqLeft
 (\mc{R}_{H}^{\alpha}, \mc{K}_{X}^{-\frac{\alpha}{2}})_{L^2(\R^{2N})} \\ 
 &= \lim_{\eps \to 0} \int_{\R^N} \left( \int_{\R^N \setminus B_{\eps}(y)} I_ {\alpha}(x-y) H(x) H(y)\left( \dive(X)(x) - (N-\alpha)\frac{(x-y)\cdot X(x)}{|x-y|^2} \right) dx \right) dy.
\end{align*}
Since $H \in Lip_{loc}(\R^N) \subset W^{1,\infty}(\supp(X))$, we have
\begin{align*}
\frac{1}{C_{N,\alpha}} (\mc{R}_{H}^{\alpha}, \mc{K}_{X}^{-\frac{\alpha}{2}})_{L^2(\R^{2N})}
=& - \lim_{\eps \to 0} \int_{\R^N} \left(\int_{\R^N \setminus B_{\eps}(y)}\frac{1}{|x-y|^{N-\alpha}} H(y) \nabla H(x) \cdot X(x) dx \right) dy + \\
& + \lim_{\eps \to 0} \frac{1}{\eps^{N-\alpha+1}} \int_{\R^N} \left(\int_{\partial B_{\eps}(y)} H(x) H(y) X(x) \cdot (x-y) d\sigma(x) \right) dy.
\\
=:& \, -\lim_{\eps \to 0} I_{\eps} +\lim_{\eps \to 0} E_{\eps}. 
\end{align*}

For $I_{\eps}$ we notice that $(I_{\alpha}*|H|)|\nabla H| |X| \in L^1(\R^N)$, thus 
$ (x,y) \mapsto I_{\alpha}(x-y) H(y) \nabla H(x) \cdot X(x) \in L^1(\R^{2N})$
and we can apply (twice) Fubini theorem; moreover $I_{\alpha}(x-\cdot) H \in L^1(\R^N)$, and we can apply Dominated Convergence theorem.
Hence we obtain
\begin{align*}
C_{N,\alpha} \lim_{\eps \to 0} I_{\eps} &= \lim_{\eps \to 0} \iint_{|x-y|>\eps} I_{\alpha}(x-y) H(y) \nabla H(x) \cdot X(x) dx dy \\
&= \lim_{\eps \to 0} \int_{\R^N} \nabla H(x)\cdot X(x) \left(\int_{\R^N \setminus B_{\eps}(x)} I_{\alpha}(x-y) H(y) dy \right) dx \\
&= \int_{\R^N} \nabla H(x) \cdot X(x) \left( \lim_{\eps \to 0}\int_{\R^N \setminus B_{\eps}(x)} I_{\alpha}(x-y) H(y) dy \right) dx \\
&= \int_{\R^N} \nabla H(x) \cdot X \, (I_{\alpha}*H).
\end{align*}

We can write $E_{\eps}$ instead as
\begin{align*}
E_{\eps} 
&= \frac{1}{2\eps^{N-\alpha+1}} \iint_{|x-y|=\eps} H(x) H(y) (X(x)-X(y)) \cdot (x-y) d\sigma(x) \times dy .
\end{align*}
If $\supp(X)\subset B_R(0)$, set $A_{R,\eps}:=\{(x,y) \in B_R(0)\times B_R(0) \mid |x-y|=\eps\}$
and observed that $H \in L^{\infty}(\R^N)$, we obtain 
$$ E_{\eps} \lesssim \frac{1}{\eps^{N-\alpha+1}} \iint_{A_{R,\eps}} |x-y|^2 d\sigma(x) \times dy = \eps^{-N+\alpha+1} m_{2N-1}(A_{R,\eps}) \lesssim \eps^{\alpha} \to 0,$$
being $\alpha>0$. This concludes the proof.
\QED

\medskip

\begin{Theorem}
\label{thm_pohoz_X}
Let $s\in (0,1)$ and $\alpha \in (0,N)$ and assume that \hyperref[(F1)]{\textnormal{(F1)}}-\hyperref[(F2)]{\textnormal{(F2)}} hold. 
Let $u \in H^s(\R^N) \cap C^{\gamma}_{loc}(\R^N) \cap Lip_{loc}(\R^N)$ for some $\gamma >2s$, be a pointwise solution of \eqref{eq_Choquard_genericaF}. 
Then $u$ satisfies the Pohozaev identity \eqref{eq:2.5}.
\end{Theorem}

\claim Proof.
We apply Proposition \ref{prop_int_kern_div} and Proposition \ref{prop_int_riesz_div} with $H=F(u)$; notice that the assumptions on $u$ and $F$ imply the needed conditions on $H$ (in particular we highlight that $f(u) \in L^{\frac{2N}{N+\alpha}}_{loc}(\R^N)$). 
Thus, for a generic $X \in C^1_c(\R^N, \R^N)$ we obtain
\begin{align*}
(\mc{G}_u^s, \mc{K}_X^s)_{L^2(\R^{2N})} & =- \big((-\Delta)^{s} u \, \nabla u,X\big)_{L^2(\R^N)} \\
&= \mu \big(u \, \nabla u,X\big)_{L^2(\R^N)} - \big((I_{\alpha}*F(u)) f(u) \, \nabla u,X\big)_{L^2(\R^N)} \\
&= \frac{\mu}{2} \big(\nabla(u^2),X\big)_{L^2(\R^N)} - \big((I_{\alpha}*F(u)) \nabla F(u),X\big)_{L^2(\R^N)} \\
&=-\frac{\mu}{2} \big(u^2,\dive(X)\big)_{L^2(\R^N)} + (\mc{R}_{F(u)}^{\alpha}, \mc{K}_{X}^{-\frac{\alpha}{2}})_{L^2(\R^{2N})}.
\end{align*}
In particular, we apply the result to 
$$X_n(x):= \varphi_n(x) \, x$$
 where $\varphi_n$ is a cut-off function with $\varphi_n \equiv 1$ in $B_n(0)$, $\supp(\varphi_n) \subset B_{n+1}(0)$, $\norm{\varphi_n}_{\infty}=1$ and $|x| |\nabla \varphi_n(x)|\leq C$ for each $x \in \R^N$ and $n \in \N$; for instance, defined such $\varphi_1$, we can set $\varphi_n:=\varphi_1(\cdot/n)$ and obtain
$$|x| |\nabla \varphi_n(x)|= |x/n| |\nabla \varphi_1 (x/n)| \leq \norm{|x| |\nabla \varphi_1(x)|}_{\infty}.$$
In particular, $x\mapsto x \varphi_n(x)$ is equi-Lipschitz. Noticed that $\dive(X_n) = N \varphi_n + \nabla \varphi_n \cdot x $
we gain
\begin{align*}
\MoveEqLeft (\mc{G}_u^s, \mc{K}_{X_n}^s)_{L^2(\R^{2N})} \\
=& \frac{C_{N,s}}{2}\int_{\R^N} \int_{\R^N} \frac{|u(x)-u(y)|^2}{|x-y|^{N+2s}}\left( \frac{N(\varphi_n(x)+\varphi_n(y))}{2} \right)dx dy -\\
&-\frac{C_{N,s}}{2}\int_{\R^N} \int_{\R^N} \frac{|u(x)-u(y)|^2}{|x-y|^{N+2s}} \left( \frac{N+2s}{2}\frac{(\varphi_n(x)\, x - \varphi_n(y) y)\cdot(x-y)}{|x-y|^2} \right) dx dy +\\
&+\frac{C_{N,s}}{2}\int_{\R^N} \int_{\R^N} \frac{|u(x)-u(y)|^2}{|x-y|^{N+2s}} \left(\frac{\nabla \varphi_n(x) \cdot x + \nabla \varphi_n(y) \cdot y}{2} \right) dx dy \\
\to& \, \frac{C_{N,s}}{2} \int_{\R^N} \int_{\R^N} \frac{|u(x)-u(y)|^2}{|x-y|^{N+2s}} \left( N- \frac{N+2s}{2}\right) = \frac{N-2s}{2} [u]_{\R^N}^2
\end{align*}
where we used $\varphi_n \to 1$, $\nabla \varphi_n \to 0$ and Dominated Convergence theorem.
Similarly
\begin{align*} 
(\mc{R}_{F(u)}^{\alpha}, \mc{K}_{X_n}^{\alpha})_{L^2(\R^{2N})} &\to \int_{\R^N} \int_{\R^N} I_{\alpha}(x-y) F(u(x)) F(u(y)) \left( N - \frac{N-\alpha}{2} \right) dx dy \\ 
&= \frac{N+\alpha}{2} \int_{\R^N} \big(I_{\alpha}*F(u)\big) F(u).
\end{align*}
and
$$\frac{\mu}{2} \big(u^2,\dive(X_n)\big)_{L^2(\R^N)} \to \mu \frac{N}{2} \norm{u}_2^2.$$
Joining the pieces, we have the claim. 
\QED

\bigskip

\medskip

\claim Proof of Theorem \ref{thm_pohoz_frac_choq}.
The theorem is a consequence of Corollary \ref{corol_reg_pohz} and Theorem \ref{thm_pohoz_X}.
\QED

\begin{Remark}
We comment the name of $\mc{K}_{X}^s$. Indeed, up to a multiplicative constant, we have, for any $\beta \in (0,1)$ and $X \in Lip_c(\R^N, \R^N)$, by \cite[equations (2.9c) and (2.11)]{CoSt} (see also \cite{{Silh0}})
\begin{align*}
\int_{\R^N} \int_{\R^N} \frac{\mc{K}_{X}^s(x,y)}{|x-y|^{N+\beta-1}} =& \int_{\R^N} \left( \int_{\R^N} \frac{\dive(X(y))}{|x-y|^{N+\beta-1}}dy\right)dx- \\
&- \frac{N+2s}{2} \int_{\R^N} \left( \int_{\R^N} \frac{(X(x)-X(y))\cdot(x-y)}{|x-y|^{N+\beta+1}}dy \right)dx \\
=& (N+\beta-1) \int_{\R^N} \dive^{\beta}(X)(x) dx - \frac{N+2s}{2} \int_{\R^N} \dive^{\beta}(X)(x) dx\\
=& (N+2\beta-2-2s) \int_{\R^N} \dive^{\beta}(X) ; 
\end{align*}
in particular
$$\int_{\R^N} \int_{\R^N} \frac{\mc{K}_{X}^s(x,y)}{|x-y|^{N+s-1}} = (N-2) \int_{\R^N} \dive^{s}(X). $$
We refer also to \cite[Chapter 3]{Dje0} where $\mc{K}_{X}^s$ is seen as the derivative of a suitable family of \emph{deformations}.
\end{Remark}


\section*{Acknowledgments}

The first and the second author are supported by PRIN 2017JPCAPN ``Qualitative and quantitative aspects of nonlinear PDEs'', and by INdAM-GNAMPA. 
The first author acknowledges financial support from PNRR MUR project PE0000023 NQSTI-National Quantum Science and Technology Institute (CUP H93C22000670006)
and PNRR MUR project CN00000013 HUB - National Centre for HPC, Big Data and Quantum Computing (CUP H93C22000450007).
The second author is supported by INdAM-GNAMPA Project “Mancanza di regolarità e spazi non lisci: studio di autofunzioni e autovalori”, codice CUP E53C23001670001.
The third author is supported in part by Grant-in-Aid for Scientific Research (22K03380, 19H00644, 18KK0073, 17H02855) of Japan Society for the Promotion of Science.


\section*{Conflict of interest}

Authors state no conflict of interest.

\end{document}